\documentclass[11pt]{article}
\usepackage[english]{babel}
\usepackage{amsmath}
\usepackage{amsthm}
\usepackage{amsfonts,amssymb,latexsym,graphicx,psfrag}
\usepackage{version}

\newtheorem{thm}{Theorem}[section]
\newtheorem{cor}[thm]{Corollary}

\newtheorem{prop}[thm]{Proposition}

\theoremstyle{definition}
\newtheorem{defn}[thm]{Definition}
\theoremstyle{remark}

\numberwithin{equation}{section}

\def \x{\times}

\def \m{\mathbb}

\def \<{\langle}
\def \>{\rangle}

\newcommand{\R}{\mathbb{R}}
\newcommand{\N}{\mathbb{N}}

\newcommand{\B}{\mathbb{B}}

\begin{document}
\title{Multifunctions of Bounded Variation}
\author{R. B. Vinter\thanks{R. B. Vinter is with Department of Electrical and Electronic Engineering, Imperial College London, Exhibition Road, London SW7 2BT, UK, {\tt\small   r.vinter@imperial.ac.uk}}}
\maketitle

\begin{abstract}
Consider control systems described by a differential equation with a control term or, more generally, by a differential inclusion with velocity set $F(t,x)$. Certain properties of  state trajectories can be derived when it is assumed that $F(t,x)$ is merely measurable w.r.t. the time variable t. But sometimes a refined analysis requires the imposition of stronger hypotheses regarding the time dependence ). Stronger forms of necessary conditions for minimizing state trajectories  can derived, for example, when $F(t,x)$ is Lipschitz continuous w.r.t. time. 
It has recently become apparent that significant addition properties of state trajectories can still be derived, when the  Lipschitz continuity hypothesis is replaced by the weaker requirement that $F(t,x)$ has bounded variation w.r.t. time. This paper introduces a new concept of multifunctions $F(t,x)$ that have bounded variation w.r.t. time near a given state trajectory, of special relevance to control. We provide an application to sensitivity analysis.

\end{abstract}

\bigskip

\noindent
{\small Keywords: Differential Inclusions, Optimal Control, Bounded Variation, Sensitivity.
\bigskip

}

\newpage
\section{Introduction}
A widely used framework for  control systems analysis is based on a description of the dynamic constraint in the form of a differential inclusion
\begin{equation}
\label{diffincl}
\dot x (t) \in F(t,x(t))\quad \mbox{a.e. }  t \in [S,T]\,,
\end{equation}
in which $F(\cdot,\cdot): [S,T]\times \R^{n} \leadsto \R^{n}$ is a given multifunction. We refer to absolutely continuous functions $x(\cdot):[S,T] \rightarrow \R^{n}$ satisfying (\ref{diffincl}) as state trajectories.
\ \\

\noindent
It is well known that the assumptions that are made regarding the $t$ dependence of $F(t,x)$ have a critical effect on the qualitative properties of the set of state trajectories and, if state trajectories minimizing a given cost function are of primary interest,  the assumptions affect the regularity properties of the value function, the nature of necessary conditions that can be derived, etc. In previous research on the distinct properties of state trajectories, depending on the different assumptions that are made about the regularity of $F(t,x)$ with respect to $t$, the attention has focused on consequences of hypothesizing:
\ \\

\noindent
(a): $t \rightarrow F(t,x)$ is measurable, or
\ \\

\noindent
(b): $t \rightarrow F(t,x)$ is Lipschitz continuous.
\ \\

\noindent
(\cite[Chap. 2]{Vinter}) definitions of measurability and Lipschitz continuity of multifunctions.) Some examples of distinct properties are as follows. 
\ \\

\noindent
(i): standard necessary conditions of optimality, in state-constrained optimal control, take a non-degenerate form, under the assumption that $F(\cdot,x)$ is Lipschitz continuous and  other assumptions, but this is no longer in general the case if $F(\cdot,x)$ is merely measurable \cite[Thm. 10.6.1]{Vinter}. 
\ \\

\noindent
(ii):  optimal state trajectories have essentially bounded derivatives under the assumption that $F(\cdot,x)$ is Lipschitz continuous and  other assumptions, but may not be essentially bounded if  $F(\cdot,x)$ is merely measurable \cite{ClarkeVinter}.
\ \\

\noindent
(iii):  the Hamitonian evaluated along an optimal state trajectory and co-state trajectory cannot contain jumps if $F(\cdot,x)$ is Lipschitz continuous, but may be discontinuous if $F(\cdot,x)$ is merely measurable \cite{Clarkebook}.
\ \\

\noindent
Other examples where there are significant differences in the implications of the two kinds of regularity hypotheses arise in the study of regularity properties of the value function for state constrained optimal control problems \cite{valuefunction}, validity of necessary conditions of optimality for free-time optimal control problems \cite{Clarkebook}, the interpretation of costate trajectories as gradients of the value function \cite{sensitivity} and in more general sensitivity analysis.
\ \\

\noindent
Are there other classes of differential inclusions $F(t,x)$, defined by their regularity w.r.t. $t$, where interesting, distinct properties are encountered?
It turns out that  multifunctions $F(t,x)$ having bounded variation w.r.t. $t$ is an example of such a class. Many properties of the set of state trajectories that are valid when $F(t,x)$ has Lipschitz $t$-dependence, but not in general when $F(t,x)$ has measurable $t$-dependence, have analogues  when $F(t,x)$ has bounded variation $t$-dependence.
%
\ \\

\noindent
How should we define `$F(t,x)$  has bounded variation $t$-dependence'? An obvious approach is to require:
\begin{equation}
\label{bounded_var}
\sup \left\{ \sum_{i=0}^{N-1} \underset{x \in X}{\sup}\,d_{H}(F(t_{i+1},x),F(t_{i},x) ) \right\}\, < \, \infty\,.
\end{equation}
Here,  $X$ is some suitably large subset of $\R^{n}$. $d_{H}(\cdot,\cdot)$ denotes the Hausdorff distance. (See \ref{HD}).  The outer supremum is taken over all possible partitions $\{t_{0}=S, \ldots, t_{N}=T\}$ of $[S,T]$. But we follow a more refined approach, for reasons that we now describe.
\ \\

\noindent
In the study of the implications of regularity assumptions regarding the $t$-dependence of $F(t,x)$, the interest usually focuses on a particular state trajectory  $\bar x(.)$ (typically a state trajectory minimizing a given cost function). We can expect that, in such situations, properties of $F(\cdot,\cdot)$ only on some neighborhood of the graph of $\bar x(.)$ would be relevant to the ensuing control systems analysis. One way to take account of the special trajectory $\bar x(\cdot)$ would be to let $X$ in (\ref{bounded_var}) be a closed set which contained  all possible values of $\bar x(.)$, i.e.
$$
\{\bar x(t) \,|\,  t \in [S,T]\} \; 
\subset\; X \quad \mbox{for all }t \in [S,T]\,.
$$
This approach involves making unnecessary assumptions about values of $F(\cdot,\cdot)$ at points far from the graph of $\bar x(\cdot)$. We therefore adopt a more refined definition of bounded variation multifunctions,  in which the inner suprema in (\ref{bounded_var}) are taken, not over $X$, but
over smaller sets (defined by a parameter $ \delta >0 $) and the outer supremum is taken over partitions $\{t_{i}\}$ of $[S,T]$  with `mesh size', which we write $\mbox{diam}\,(\{t_{i}\})$ (see (\ref{mesh}) below),  not greater than $\epsilon >0$. Accordingly, we say that $t\rightarrow F(t,\cdot)$ has bounded variation along $\bar x(\cdot)$ if, for some $\delta >0$ and $\epsilon >0$ we have
\begin{equation}
\label{supremum2}
\sup \left\{ \sum_{i=0}^{N-1} \underset{x \in \bar x(t)+\delta \B, t \in [t_{i},t_{i+1}]}{\sup}\,d_{H}(F(t_{i+1},x),F(t_{i},x) ) \,|\,  \mbox{diam}( \{t_{i}\}) \leq \epsilon \right\}\, < \, \infty\,.
\end{equation}
We add another refinement; that is to consider multifunctions $F(t,x,a)$, whose argument includes an additional variable $a$ that ranges over a given subset $A$ of a finite dimensional linear space. Including the parameter $a$ provides useful flexibility for certain applications \cite{Palladino}.
\ \\

\noindent
{\it The purpose of this paper is to bring together and prove  properties (relevant to control system analysis) of a multifunction that has bounded variation along some given state trajectory $\bar x(\cdot)$, and of the associated cummulative variation function.} These include one-sided continuity properties of such multifunctions and the effects on the cummulative variation function of changes to the multifunction.  In the case of a function $m(\cdot,x)$ of bounded variation along $\bar x(\cdot)$ (a function can be regarded as a special case of a multifunction),  it is shown that there is an associated signed Borel measure. Finally, we show how this theory can be used to obtain new sensitivity formulae describing how the output of a control system is affected by a small time delay in the implementation of a control.  
\ \\

\noindent
The analysis in this paper generalizes some aspects of the classical theory of functions of a scalar variable having bounded variation,  to allow for several independent variables, when the bounded variation property pertains only to one of the variables, and when multifunctions replace functions. There is extensive recent work, treating the properties of bounded variation functions with several independent variables, for which the monograph \cite{Ambrosio} is a comprehensive source of references. The motivation arises from a desire to investigate regularity properties of minimizers of variational problems in several independent variables and of solutions to Hamilton Jacobi equations arising in optimal control (see, for example, \cite{Tonon}). Multi-functions $F(t)$ of a single variable $t$ (no $x$-dependence) possessing a one-sided bounded variation property have been investigated by Moreau \cite{Moreau}, in connection with sweeping processes. But the study initiated in this paper of multifunctions that are $x$-dependent and have bounded variation `near'  a given state trajectory $\bar x(.)$ is a new development.
\ \\

\noindent
{\it Notation:} 
For vectors $x\in \m R^{n}$, $|x|$ denotes the Euclidean length. $\B$ denotes the 
closed
unit ball in $\m R^{n}$. 
Given a multifunction $\Gamma(\cdot): \m R^{n} \leadsto \m R^{k}$, the graph of $\Gamma(\cdot)$, written Gr$\,\{\Gamma(\cdot)\}$, is the set $\{(x,v)\in \m R^n\x\m R^k\,|\,v\in \Gamma(x)\}$. Give a set $A \subset \m R^{n}$ and a point $x \in \m R^{n}$, we denote by $d_{A}(x)$ the Euclidean distance of a point $x\in \m R^{n}$ from $A$:
$$
d_{A}(x)\;:=\; \inf\{|x-y| \,|\, y \in A \}\;.
$$ 
co\,$A$ denotes the convex hull of a set $A \in \R^{n}$. Given an interval $I$,  we write $\chi_{I}(t)$ for the indicator function of $I$, taking values $1$ and $0$ when $t\in I$ and $t \notin I$, respectively. For numbers $a$ an $b$, $a \vee b := \max\{a,b\}$ and
$a \wedge b := \min\{a,b\}$. Given two nonempty sets $A, A \in \R^{k}$, their Hausdorff distance is
\begin{equation}
d_{H(A,B)}\,:=\, \inf\{d_{A}(x)| x\in B\} \vee \inf\{d_{B}(x)| x\in A\}\,.
\label{HD}
\end{equation}
A function $r:[S,T] \rightarrow \R$ of bounded variation on the interval $[S,T]$ has a left limit, written $r(t^{-})$, at every point $t \in (S,T]$ and  a right limit, written $r(t^{+})$, at every point $t \in [S,T)$. We say $r(.)$ is normalized if it is right continuous on $(S,T)$. 
\ \\

\noindent
We denote by $NBV^{+}[S,T]$ the space of increasing, real-valued functions $\mu(.)$ on $[S,T]$  of bounded variation, vanishing at the point $S$ and right continuous on $(S,T)$. The total variation of a function $\mu(\cdot)\in NBV^{+}[S,T]$ is written $||\mu(\cdot)||_{\mbox{TV}}$. As is well known, each point $\mu(\cdot)\in NBV^{+}[S,T]$ defines a unique Borel measure on $[S,T]$. This associated measure is also denoted $\mu(\cdot)$. The space of continuous functions $x(\cdot):[S,T]\rightarrow \R^{n}$ with supremum norm is written $C([S,T];\R^{n})$ and we denote by $C^{*}([S,T];\R^{n})$ its topological dual space.
\ \\

\noindent
A {\it modulus of continuity} is a function $\theta(.): [0,\infty)\rightarrow [0, \infty)$ such that $\lim_{s \downarrow 0} \theta(s)=0$.
\ \\

\noindent
%
Take a lower semicontinuous function $f(\cdot):\m R^k\to \m R\cup\{+\infty\}$ and a point $\bar x \in \mbox{dom}\, f(\cdot) \,:=\, \{x\in \R^{k}\,|\, f(x) < + \infty\}$. 
The {\it subdifferential of $f(\cdot)$ at $\bar x$}, 
denoted $\partial f(\bar x)$ is the set:
\begin{eqnarray*}
\lefteqn{\partial f(\bar{x}) \;:=\;  \Big\{ \xi\,|\, \exists\; \xi_{i} \rightarrow \xi \mbox{ and }
x_{i} \stackrel{\mbox{{\tiny {\rm dom}}}\,f(\cdot)}{\longrightarrow} \bar{x} \mbox{ such that}
}
\\[3mm]
&& \limsup_{x \rightarrow  x_{i}} \frac{ \xi_{i} \cdot (x-x_{i})- \varphi(x)+\varphi(x_{i})}{|x-x_{i}|}
\, \leq  0 \, \mbox{ for all } i \in \N \Big\} \;. 
\end{eqnarray*}
Here, `$x_{i} \stackrel{\mbox{dom}\,f(\cdot)}{\rightarrow} \bar x$' means that all elements in the convergent sequence $\{x_{i}\}$ lie in $\mbox{dom}\,f(\cdot)$.
For further information about subdifferentials, and related constructs in nonsmooth analysis, see \cite{clarke}, \cite{rockafellar} and \cite{Vinter}.
\section{Multifunctions of Bounded Variation}

%
%
%
Take a bounded interval $[S,T]$, a compact set $A\subset \R^{k},$ a multifunction $F(.,.,.): [S,T]\times \R^{n} \times A\leadsto \R^{n}$ and a continuous function $\bar{x}(.): [S,T]\rightarrow \R^{n}$. Generic elements in the domain of $F(\cdot,\cdot,\cdot)$ are denoted by $(t,x,a)$.  
\ \\

\noindent
In this section we define a concept that makes precise the statement `$F(t,x,a)$ has bounded variation with respect to the $t$ variable, along  $\bar{x}(.)$, uniformly with respect to $a\in A$'. If $F(t,x,a)$ is independent of $(x,a)$ and single valued, i.e. $F(t,x,a)=\{f(t)\}$ for some function $f(\cdot):[S,T]\rightarrow \R^{n}$, this concept reduces to the standard notion  `$f(\cdot)$ has bounded variation'.
\ \\

\noindent
For any $t \in [S,T]$, $\delta >0$ and partition 
$
{\cal T}\,=\, \{t_{0}=S,t_{1},\ldots,t_{N-1}, t_{N}=t\}
$
of $[S,t]$, 
define $I^{\delta}({\cal T})\in \R^{+} \cup \{+ \infty\}$ to be 
$$
I^{\delta}({\cal T})\,:=\, 
\sum_{i=0}^{N-1}\sup \left\{ d_{H}(F(t_{i+1},x,a) ,F(t_{i},x,a) )\;|\; x \in \bar x([t_{i},t_{i+1}])+\delta B, a\in A  \right\}\;.
$$
Here, $\bar x([t_{i},t_{i+1}])$ denotes the set $\{\bar x(t)\,|\, t \in  [t_{i},t_{i+1}]\}$.
\ \\

\noindent
Take any $\epsilon >0$.  Let $\eta_{\epsilon}^{\delta}(\cdot): [S,T] \rightarrow \R^{+} \cup \{+ \infty\}$ be the function defined as follows: $\eta_{\epsilon}^{\delta}(S)=0$ and, for $t \in (S,T]$,
\begin{equation}
\label{eps-delta}
\eta_{\epsilon}^{\delta}(t)\;=\; \sup \left\{ I^{\delta}({\cal T})\,|\, {\cal T} \mbox{ is a partition of } [S,t] \mbox{ s.t. } \mbox{diam}({\cal T})\leq \epsilon 
 \right\}\,,
\end{equation}
in which
\begin{equation}
\label{mesh}
\mbox{diam} ( {\cal T} )\,:=\, \sup \{ t_{i+1}-t_{i} \,|\, i=0,\ldots,N-1 \}\,.
\end{equation}
It is clear that, for any $t \in [S,T]$, $\delta >0$, $\delta' >0$, $\epsilon >0$, $\epsilon' >0$,  
\begin{equation}
\label{elementary}
\delta' \leq \delta \mbox{ and }  \epsilon' \leq \epsilon \implies 
0\leq \eta_{\epsilon'}^{\delta'}(t) \leq  \eta_{\epsilon}^{\delta}(t)\;.
\end{equation}
 (This relation is valid even when $\eta_{\epsilon}^{\delta}(t)=+  \infty$, according to the rule  `$+\infty \leq +\infty$'.) We may therefore define the functions $\eta^{\delta}(.), \eta(.): [S,T] \rightarrow \R^{+} \cup \{+ \infty\}$ to be
\begin{eqnarray}
\label{1.1z}
&&\eta^{\delta}(t)\,:=\, \lim_{\epsilon \downarrow 0} \eta^{\delta}_{\epsilon}(t) \;\mbox{ for } t \in [S,T]\,
\\
\label{1.2z}
&&\eta(t)\,:=\, \lim_{\delta \downarrow 0} \eta^{\delta}(t) \;\mbox{ for } t \in [S,T]\,.
\end{eqnarray}

\noindent
\begin{defn} 
\label{definition}
Take a  set $A \subset \R^{k}$, a multifunction $F(\cdot,\cdot,\cdot): [S,T]\times \R^{n} \times A\leadsto \R^{n}$ and a  function $\bar{x}(\cdot): [S,T]\rightarrow \R^{n}$.
We say that {\it $F(\cdot,x,a)$ has bounded variation along $\bar{x}(.)$ uniformly over  $A$}, if the function $\eta(\cdot)$ given by (\ref{1.2z}) satisfies
$
\eta(T)\,< \, +\infty\;.
$
\end{defn}

\noindent
If $F(\cdot,x,a)$ has bounded variation along $\bar{x}(\cdot)$ uniformly over  $A$,   the function $\eta(\cdot)$ defined by (\ref{1.2z}) is called  the {\it cummulative variation function} of $ F(\cdot,x,a)$ along $\bar{x}(\cdot)$, uniformly over $A$. We also refer to $\eta^{\delta}_{\epsilon}(\cdot)$ and $\eta^{\delta}(\cdot)$, defined by (\ref{eps-delta}) and (\ref{1.1z}) as the $(\delta,\epsilon) $-perturbed cummulative variation function and $\epsilon$-perturbed cummulative variation function respectively.
\ \\

\noindent
In what follows we will adhere to the following notational convention: if $\eta(\cdot)$ is a given cummulative variation function (for some $F(\cdot,\cdot,\cdot)$ and $\bar x(\cdot)$) then, for any $\delta >0$ and $\epsilon >0$, $\eta^{\delta}(\cdot)$ and $\eta^{\delta}_{\epsilon}(\cdot)$  denotes the $\delta$-perturbed and $\delta,\epsilon$-perturbed cummulative variation functions associated with $\eta(\cdot)$, according to (\ref{eps-delta}) and (\ref{1.1z}).
\ \\

\noindent
If $F(t,x,a)$ does not depend on $a$, we omit mention of the qualifier `uniformly over $A$'. A function $ L(\cdot,x,a)$ is said to have bounded variation along $\bar{x}(.)$ uniformly over  $A$, if the associated multifunction $ \{L(\cdot,x,a) \}$ has this property.
%
\ \\

\noindent
Assume that $ F(\cdot,x,a)$ has bounded variation along $\bar x(\cdot)$ uniformly over $A$.   Then there exist $\bar{ \delta}>0$ and $\bar{\epsilon}>0$ for which $\eta^{\bar \delta}_{\bar \epsilon}(T) < +\infty$. We list the following elementary properties of the accumulative variation functions (`elementary', in the sense that they are simple consequences of the definitions): for any $\delta \in (0,\bar \delta]$ and $\epsilon \in (0,\bar \epsilon]$,

\begin{itemize}
\item[(a):] $t \rightarrow \eta^{\delta}_{\epsilon }(t)$, $t \rightarrow \eta^{\delta}(t)$ and $t \rightarrow \eta(t)$ are  increasing, finite valued functions,
\item[(b):]  $\eta^{\delta}_{\epsilon }(t) \geq \eta^{\delta}(t) \geq \eta(t)$ for all $t \in [S,T]$
\item[and]
\item[(c):] given any $[s,t]\subset [S,T]$ such that $t-s \leq \epsilon$,
\begin{equation}
\label{basic}
d_{H}(F(t,y,a),F(s,y,a)) \,\leq \,
\eta^{\delta}_{\epsilon}(t)-\eta^{\delta}_{\epsilon}(s), 
\end{equation}
$\mbox{ for all }y \in \bar{x}([s,t])+ \delta \B\,\mbox{ and } \, a\in A$.
\end{itemize}
%
%

\noindent
{\bf Example.} An important potential role of the preceding constructs will be to derive regularity  properties of value functions, minimizing state trajectories and other functions associated with an optimal control problem, in which the dynamic constraint is a differential inclusion $\dot x \in F(t,x)$, when $F(t,x)$ has bounded variation with respect to the $t$ variable, `near' a given state trajectory $\bar x(\cdot)$. Regularity properties are typically related to the cummulative variation function $\eta(.)$ of  $F(\cdot,x)$. The more precise is the information about the cummulative variation the more informative is the corresponding regularity property that can be derived.
 This is the main reason why we have adopted the refined definition, Def. \ref{definition},  for the formulation of the `bounded variation' hypothesis, in place of a simpler one based on the condition (\ref{bounded_var}), for some closed subset $X $ that strictly contains the range of $\bar x(\cdot)$ in its interior.
The purpose of this example is to show that using the `refined'   definition can provide a more informative cummulative variation function.
\ \\

\noindent
Consider the function $f(\cdot,\cdot):[0,1]\times \R \rightarrow \R$ and the function $\bar x(\cdot):[0,1] \rightarrow \R$:
$$
f(t,x)=t\, x \mbox{ and } \bar x(t) =t \quad \mbox{for $(t,x) \in [0,1]\times \R$.}
$$

Take $X= \mbox{range}\{\bar x(\cdot)\}  = [0, 1]$. The cummulative variation function $\eta_{\mbox{simple}}(\cdot)$ of $ f(\cdot,x)$ related to condition (\ref{bounded_var}) and defined by 
$$
\eta_{\mbox{simple}}(t) \,=\, \sup \left\{ \sum_{i=0}^{N-1} \underset{x \in X}{\sup}\,d_{H}(F(t_{i+1},x),F(t_{i},x) ) \right\}
$$  
in which the outer supremum is taken over all possible partitions of $[0,t]$, is easily calculated to be:  
$$
\eta_{\mbox{simple}}(t)\,=\, t\,.
$$
Also, the cummulative variation of $f(\cdot,x)$ along $\bar x(.)$ according to Def. \ref{simple}, is
$$
\eta(t)\,=\, \frac{1}{2} t^{2}\,.
$$
Notice that for any (nontrivial) subinterval $[s,t] \subset [0,1]$, $t >s$, 
\begin{eqnarray*}
&&(\eta_{\mbox{simple}}(t)- \eta_{\mbox{simple}}(s) ) -
(\eta(t)- \eta(s) )
\;=
\\
&& \hspace{1.0 in}
(t-s) - \frac{1}{2}(t^{2}-s^{2}) = (t-s)(1- \frac{1}{2}(t+s)) \,(\,>\,0),  
\end{eqnarray*}
from which it can be deduced that the Borel measure induced by $\eta(\cdot)$ strictly minorizes that induced by $\eta_{\mbox{strict}}(\cdot)$ in the sense
$$
\int_{D} d\eta_{\mbox{simple}}(t) - \int_{D} d\eta(t)
 >0\, 
$$
for any Borel subset $D \subset [0,1]$ having nonempty interior. This demonstrates the greater precision that can be achieved in regularity analysis, by using the more refined definition.
\section{Continuity Properties}
As is well known, an $\R^{n}$-valued function of bounded variation on a finite interval may be discontinuous, but it has everywhere left and right limits and it has at most a countable number points of discontinuity. A multifunction having bounded variation along  a given continuous trajectory uniformly over a given set has similar properties, as described in the following proposition. We invoke the hypotheses:
\begin{itemize}
\item[(C1)]
$F(t,x,a)$ is closed and non-empty for each $(t,x,a)\in [S,T]\times \R^{n}\times A$, $F(\cdot,x,a)$ is measurable for each $(x,a) \in \R^{n}\times A$ and there exists $c>0$ 
such that
\begin{equation}
\label{limsup}
F(t,x,a)\,\subset\,c\B \mbox{ for all $x \in \bar{x}(t)+  \bar{\delta}\B$, $t \in [S,T]$, $a\in A$. }
\end{equation}
\item[(C2)]
There exists a modulus of continuity $\gamma(.): \R^{+} \rightarrow \R^{+}$ such that
\begin{equation}
\label{lcontinuous}
F(t,x,a)\,\subset\, F(t,x',a') + \gamma( |x-x'|+|a-a'|) \B
\end{equation}
\mbox{ for all $x,x' \in \bar{x}(t)+  \bar{\delta}\B$, $t \in [S,T]$ and $a,a'\in A$. }
\end{itemize}

\begin{prop}\label{prop1}
Take  a compact set  $A \subset \R^{k}$, a continuous function $\bar{x}(\cdot):[S,T]\rightarrow \R^{n}$ and a  multifunction $F(\cdot,\cdot,\cdot): [S,T]\times \R^{n}\times A \leadsto \R^{n}$. Suppose that $F(\cdot,x,a)$ has bounded variation along  $\bar{x}(\cdot)$ uniformly over $A$. Take $\bar \delta >0$ such that $\eta^{\bar \delta}(T)< + \infty$. Assume (C1) and (C2). Take any $\delta \in (0,\bar \delta)$. Then
\begin{itemize}
\item[(a):]For any $\bar s \in [S,T) $ and $\bar t \in (S,T] $ , the one-sided, set-valued limits
$$
F(\bar{s}^{+},x,a)\,:=\,\underset{s\downarrow \bar s}{\lim} \,F(s,x,a)
 \,, \quad
F(\bar{t}^{-},y,a)\,:=\,\underset{t\uparrow \bar t}{\lim}\, F(t,y,a)
$$ 
exist for every $x \in \bar{x}(\bar{s})+ \delta \B$, $y \in \bar{x}(\bar{t})+ \delta \B$ and $a\in A$.  
\item[(b):]
For any $\bar s \in [S,T) $ and $\bar t \in (S,T] $ 
$$
\underset{s\downarrow \bar s}{\lim} \;
\underset{\underset{a\in A}{x \in \bar{x}(\bar s)+\delta B,}}{\sup}\,
d_{H}(F(\bar{s}^{+},x,a), F(s,x,a))\,=\,0
$$ 
and 
$$
\underset{t\uparrow \bar t}{\lim} \;
\underset{\underset{a\in A}{x \in \bar{x}(\bar t)+\delta B,}}{\sup}\,
d_{H}(F(\bar{t}^{-},x,a), F(t,x,a))\,=\,0
$$ 
\item[(c):]
There exists a countable set ${\cal A}$ such that, for every $t \in (S,T)\backslash {\cal A}$, 
$$
\underset{t' \rightarrow t}{\lim} \;
\underset{\underset{a\in A}{x \in \bar{x}(t)+ \bar{\delta} B}}{\sup}\,
d_{H}(F(t',x,a), F(t,x,a))\,=\,0\,.
$$ 
\end{itemize}
\end{prop}

\noindent
{\bf Proof.}

\noindent
(a):
We prove only the first assertion. The proof of the second assertion is similar. Choose any $\bar s \in [S,T)$. Take $\epsilon >0$ such that $\eta^{\bar \delta}_{\epsilon}(T) < + \infty$. Fix $\delta \in (0,\bar{\delta})$. Take any $x \in \bar{x}(\bar{s})+ \delta \B$, $a\in A$ and
$$
v \in \underset{s\downarrow \bar{s}}{\mbox{lim sup}}\, F(s,x,a)\;.
$$
By definition of `lim sup', there exists $s_{i} \downarrow \bar{s}$ and $v_{i}\rightarrow  v$ such that
$$
v_{i}\in F(s_{i},x,a)    \mbox{ for all }i \;\mbox{ and }\;    v_{i}\rightarrow v \mbox{ as }  i\rightarrow \infty\,.
$$ 
The assertion (a) will follow if we can show that, also,
\begin{equation}
\label{lim inf}
v \in \underset{s\downarrow \bar{s}}{\mbox{lim inf}}\, F(s,x,a)\;,
\end{equation}
i.e. the `lim sup' and `lim inf' coincide, in which case the limit exists. 
To show (\ref{lim inf}) we  take an arbitrary sequence $t_{j}\downarrow \bar{s}$. 
Since $\bar x(\cdot)$ is continuous and $x \in \bar x(\bar s)+ \delta \B$, we can arrange, by eliminating elements in the sequence $\{(s_{i},v_{i})\}$,  that, for every $j$,
$ \bar{s }\leq s_{j} < t_{j},\, t_{j}-\bar s\leq \epsilon \mbox{ and } x \in \bar{x}(t)+ \bar \delta \B\;\mbox{ for all } t \in [\bar{s},{t}_{j}]
$, $j =1,2,\ldots$
But then, since $t_{j}-s_{j}\leq \epsilon$ and by property (\ref{basic}) of the $(\delta,\epsilon)$-perturbed cummulative variation function,
$$
d_{H}(F(t_{j},x,a), F(s_{j},x,a))\,\leq\, \eta_{\epsilon}^{\bar{\delta}}(t_{j})- \eta_{\epsilon}^{\bar{\delta}}(s_{j})\;.
$$
This means that, for each $j$, there exists $w_{j} \in F(t_{j},x,a)$ and 
$$
|v_{j}-w_{j}|\,\leq\, \eta_{\epsilon}^{\bar{\delta}}(t_{j})- \eta_{\epsilon}^{\bar{\delta}}(s_{j})\;.
$$
We know however that, since $\eta_{\epsilon}^{\bar \delta}(.)$ is a finite valued, monotone function, it has a right limit $\eta_{\epsilon}^{\bar \delta}(\bar{s}^{+})$ at $\bar{s}$. Hence
$$
\lim_{j\rightarrow \infty} |v_{j}-w_{j}|\,\leq\, 
\lim_{j \rightarrow \infty}\left( \eta_{\epsilon}^{\bar{\delta}}(t_{j})- \eta_{\epsilon}^{\bar{\delta}}(s_{j})\right) \,\leq\,
\eta_{\epsilon}^{\bar{\delta}}(\bar{s}^{+})- \eta_{\epsilon}^{\bar{\delta}}(\bar{s}^{+})=0\;.
$$
It follows that $v_{j}-w_{j}\rightarrow 0$. But then $v= \lim_{j}v_{j}= \lim_{j} w_{j}$. Since $t_{j}\downarrow \bar s$ was an arbitrary sequence, we conclude (\ref{lim inf}). We have confirmed (a).
\ \\

\noindent
(b) These assertions may be deduced from (a), together with the compactness of the set $A$ and of closed $\delta$ balls about $\bar{x}(\bar{s})$ and $\bar{x}(\bar{t})$, and with the assumed continuity properties of $(x,a)\rightarrow F(t,x,a)$.
\ \\

\noindent
(c) Let ${\cal A}$ be the empty or countable subset of $(S,T)$ comprising points at which the finite-valued, monotone function $\eta_{ \epsilon}^{\bar \delta}(\cdot)$ is discontinuous. Fix a point $t \in (S,T)\backslash {\cal A}$.
Take any $\rho >0$. Since $\eta_{\epsilon}^{\bar \delta}(\R^{n}\times A)$ is continuous at $t$, we may choose $\sigma >0$ such that 
$$
\eta_{\epsilon}^{\bar \delta}(t+\sigma)-\eta_{\epsilon}^{\bar \delta}(t-\sigma) \,\leq\, \rho\;.
$$
So by (\ref{basic}), for any $t' \in [S,T]$ such that $|t'-t| \leq \sigma\wedge \epsilon$ we have 
\begin{eqnarray*}
\sup\{d_{H}(F(t',x,a),F(t,x,a))\,|\,x\in \bar{x}(t)+ \bar \delta \B,\, a\in A \}&\leq&  \eta_{\epsilon}^{\bar \delta}(t' \vee t)-\eta_{\epsilon}^{\bar \delta}(t' \wedge t)
\\
&\leq& \eta_{\epsilon}^{\bar \delta}(t+\sigma)-\eta_{\epsilon}^{\bar \delta}(t-\sigma)\leq\rho\;.
\end{eqnarray*}
The continuity properties of $F(.,x,a)$ at $t$ have been confirmed. $\square$
\ \\

\noindent
The following proposition provides information about how the cummulative variation function of a multifunction, and its $\delta$-perturbation, are affected by changes of $\delta$ and the parameter space for $a$. 
\begin{prop}
\label{simple}
Take compact sets $A_{1}, A \subset \R^{k}$ such that $A_{1} \subset A$, a continuous function $\bar{x}(.)$ and a multifunction $F(\cdot,\cdot,\cdot): [S,T]\times \R^{n}\times \R^{k}\leadsto \R^{n}$. Suppose that $ F(\cdot,x,a)$ has bounded variation along $\bar{x}(.)$ uniformly over $A$. Write the cummulative variation function $\eta_{A}(\cdot)$.  Let $\bar \delta >0$ be such that $\eta^{\bar \delta}_{A}(T) < \infty$. Suppose that hypotheses (C1) and (C2) of Prop. \ref{prop1} are satisfied.
\vspace{0.05 in}

\noindent
Then $F(\cdot,x,a)$ has bounded variation along $\bar{x}(.)$ uniformly over $A_{1}$. (Write the cummulative variation function $\eta_{A_{1}}(\cdot)$.)
%
Furthermore, for any $[s,t] \subset [S,T]$ and $\delta >0$ and $\delta' >0$ such that $0  <   \delta' \leq \delta \leq \bar \delta$,

\begin{equation}
\label{1star}
%
\eta_{A_{1}}(t)-\eta_{A_{1}}(s)\;\leq\;\eta^{\delta'}_{A_{1}}(t)-\eta^{\delta'}_{A_{1}}(s)\;\leq\; \eta^{\delta}_{A}(t)-\eta^{\delta}_{A}(s)\;.
\end{equation}
\end{prop}

\noindent
{\bf Proof.} 
 The facts that  $ F(\cdot,x,a)$ has bounded variation along $\bar{x}(\cdot)$ uniformly over $A_{1}$ and  $\eta^{\bar \delta}_{A_{1}}(T) < \infty$, follow immediately from the definitions, since $A_{1} \subset A$.
\ \\

\noindent
Take any $[s,t] \subset [S,T]$ and $0< \delta'\leq \delta \leq\bar \delta$, such that $\eta^{\bar \delta}_{A}(T) < \infty $. Choose $\epsilon >0$ such that $\eta_{A\, \epsilon}^{\bar \delta}(T) < \infty$. 
\ \\

\noindent
Concerning assertion (\ref{1star}), we need prove only  the right-side inequality
%
%
$\eta^{\delta'}_{A_{1}}(t)-\eta^{\delta'}_{A_{1}}(s)\;\leq\; \eta^{\delta}_{A}(t)-\eta^{\delta}_{A}(s)$, since the left-side inequality of (\ref{1star}) follows immediately by passing to the limit as $\delta' \downarrow 0$. 
\ \\

\noindent
Presently, we shall make use of the following fact: take $[\sigma_{1},\sigma_{2}] \subset [S,T]$ such that $|\sigma_{1}-\sigma_{2}| \leq \epsilon$. Then, for any  $\sigma \in [\sigma_{1},\sigma_{2}]$ and 
$x' \in \bar{x}([\sigma_{1},\sigma_{2}])+ \delta' \B$,
\begin{eqnarray}
\nonumber
&& \sup\{ d_{H}(F(\sigma_{1},x',a),F(\sigma_{2},x',a))| a \in A\}
\\
\nonumber
&& \, \leq   \sup\{ d_{H}(F(\sigma_{1},x,a),F(\sigma,x,a))| x\in \bar{x}([\sigma_{1},\sigma])+ \delta' \B, a\in A \}
\\
&& \; +  \sup\{ d_{H}(F(\sigma,x,a),F(\sigma_{2},x,a))| x\in \bar{x}([\sigma,\sigma_{2}])+ \delta' \B, a\in A\} + 2 \gamma(c\epsilon)
\label{measure}
\end{eqnarray}
in which $c$  and $\gamma(\cdot)$ are as in hypotheses (C1) and (C2). Define 
$$
G^{\delta'}_{A_{1}\epsilon}(\sigma_{1},\sigma_{2})\,:=\,
\sup\, 
\sum_{i=0}^{N}  \sup \{
d_{H}(F(t_{i+1},x,a), F(t_{i},x,a)) |    x \in \bar x([t_{i},t_{i+1}]) + \delta' \B, a \in A_{1} \}
\,.
$$
(The outer supremum is taken over partitions $\{t_{i}\}$ of $[\sigma_{1},\sigma_{2}]$ with  diam $(\{ t_{i}\}) \leq\epsilon$.)
\vspace{0.05 in}

\noindent
In view of (\ref{measure}) and by consideration of arbitrary partitions of $[S,t]$ having diameter at most $\epsilon$, and their modifications to include an extra grid point $s$, we can show
\begin{equation}
\label{first_ineq}
\eta^{\delta'}_{A_{1}\,\epsilon} (t)\,\leq\, 
\eta^{\delta'}_{A_{1}\,\epsilon}(s)+ 2 \, \gamma (c \,\epsilon) + G^{\delta'}_{A_{1} \epsilon}(s,t)\,,
\end{equation}
 From the definition of $\eta^{\delta}_{A\, \epsilon}(.)$ 
\begin{equation}
\label{second_ineq}
\eta^{\delta}_{A\, \epsilon}(t) \,\geq\,  \eta^{\delta}_{A\,\epsilon}(s)+G^{\delta}_{A\epsilon}(s,t)
\;.
\end{equation}
Since $A_{1} \subset A$ and $\delta' \leq \delta$, $G^{\delta'}_{A_{1}\epsilon}(s,t)\leq G^{\delta}_{A\epsilon}(s,t)$.   Using this relation, combining (\ref{first_ineq}) and (\ref{second_ineq}) and passing to the limit as $\epsilon \downarrow 0$,
yields
$$
\eta^{\delta'}_{A_{1}}(t) -\eta^{\delta'}_{A_{1}}(s)
\,\leq\, \eta^{\delta}_{A}(t) -\eta^{\delta}_{A}(s)\,.\qquad \square
$$ 

\noindent
The next proposition  relates the cummulative variation function of the multifunction $F(\cdot,x,a)$ to that of the derived multifunction $\tilde F(\cdot,x,a)$, obtained by replacing the end-point values by  left and right limits.
\begin{prop}
\label{prop2.3}
Take  a compact set  $A \subset \R^{k}$, a continuous function $\bar{x}(\cdot):[S,T]\rightarrow \R^{n}$ and a  multifunction $F(\cdot,\cdot,\cdot): [S,T]\times \R^{n}\times A \rightarrow \R^{n}$. Suppose that $F(\cdot,x,a)$ has bounded variation along  $\bar{x}(\cdot)$ uniformly over $A$.  Denote by $\eta(\cdot)$  the cummulative variation function. Take $\bar \delta >0$ such that $\eta^{\bar \delta}(T) < \infty$.
Assume that hypotheses (C1) and (C2).
Let $\tilde{F}(\cdot,\cdot,\cdot): [S,T]\times \R^{n}\times A\rightarrow \R^{n}$ be a multifunction such that, for $(t,x,a)\in [S,T] \times \R^{n}\times A$,
\begin{equation}
\label{modified}
\tilde{F}(t,x,a)\,=\,
\left\{
\begin{array}{ll}
F(S^{+},x,a)& \mbox{if } t=S \mbox{ and } |x-\bar{x}(S)| < \bar \delta
\\
F(T^{-},x,a)& \mbox{if } t=T \mbox{ and }  |x-\bar{x}(T)| < \bar \delta
\\
F(t,x,a)& \mbox{otherwise . } 
\end{array}
\right.
\end{equation}
(The limit sets $F(S^{+},x,a)$ and $F(T^{-},x,a))$ exist, by the preceding proposition.) 
\ \\

\noindent
Then  $\tilde{F}(\cdot,\cdot,\cdot)$ has bounded variation along $\bar{x}(\cdot)$ uniformly over $A$. Write $\tilde \eta(\cdot)$ 
for its  
cummulative variation function.
\ \\

\noindent
Take any  $\delta \in (0, \bar  \delta)$. Then $\tilde \eta^{\delta}(.)$ 
is right continuous at $S$ and left continuous at $T$, i.e.
\begin{equation}
\label{z1}
\tilde \eta^{\delta}(S)=\lim_{s\downarrow S}\tilde \eta^{\delta}(s)\quad \mbox{and} \quad \tilde \eta^{\delta}(T)=\lim_{t \uparrow T} \tilde\eta^{\delta}(t)\,.
\end{equation}
Furthermore,
the $\delta$-perturbed cummulative variation functions of $F(.,.,.)$ and $\tilde F(.,.,.)$ 
are related as follows:
\begin{equation}
\label{z2}
\tilde \eta^{\delta}(t) \,=\,
\left\{
\begin{array}{ll}
\eta^{\delta}(t) - \underset{x \in \bar x(S)+ \delta \B,\, a\in A} {\sup}\, d_{H}(F(S,x,a), F(S^{+},x,a) )
& \mbox{for } t \in (S,T)
\\
\eta^{\delta}(T) - \underset{x \in \bar x(S)+ \delta \B,\, a\in A} {\sup}\, d_{H}(F(S,x,a), F(S^{+},x,a) )
&
\\
\hspace{0.8 in}-
\underset{x \in \bar x(T)+ \delta \B,\, a\in A} {\sup}\, d_{H}(F(T,x,a), F(T^{-}, x,a) )
& \mbox{for } t =T\,.
\end{array}
\right.
\end{equation}
Relation (\ref{z2}) implies, in particular, that
\begin{equation}
\label{z3}
\tilde \eta^{\delta}(t)-\tilde \eta^{\delta}(s) = 
\eta^{\delta}(t)-\eta^{\delta}(s)\, \; \mbox{for any }[s,t] \subset
 (S,T)\,.
\end{equation}
The relations (\ref{z1})-(\ref{z3}) remain valid when  $\delta=0$, under the  interpretation `\,$\eta^{0}(\cdot)= \eta(\cdot)$ and 
$\tilde \eta^{0}(\cdot)= \tilde \eta(\cdot)$'.
\end{prop}
\vspace{0.05 in} 
%

\noindent
The  proof of Prop. \ref{prop2.3} is given  in the Appendix.

\section{The Partial Variation Measure of a Function of a Scalar and a Vector Variable}
In this section we examine in more detail the properties of a function $m(\cdot,\cdot):[S,T] \times \R^{n} \rightarrow \R^{r}$ that has bounded variation with respect to the first variable, along a given trajectory $\bar x(\cdot): [S,T] \rightarrow \R^{n}$. We restrict attention to  a special case of the multifunctions earlier considered, in which the multifunction is point valued (i.e. $m(\cdot,\cdot)$ is a function), and no longer depends on the parameter $a$. 
\ \\

\noindent
The motivation for these investigation is a desire to make sense of integrals arising in sensitivity analysis,  of the form
\begin{equation}
\label{partial_variation}
\int_{S}^{T} p^{T}(t) \frac{\partial m}{\partial t}\, (t,\bar{x}(t))dt\,,
\end{equation}
in circumstances when $m(t,x)$ has bounded variation with respect to the first variable, but fails to be continuously differentiable with respect to this variable. Here, $p(\cdot)$ is a given continuous function. Notice that, if $m(\cdot,\cdot)$ is a continuously differentiable function, the integral can be 
written as
\begin{equation}
\label{partial_variation1}
\int_{S}^{T} p^{T}(t) d \mu (t)\,,
\end{equation}
where $\mu(\cdot)$ is the (signed) Borel measure on $[S,T]$ defined by 
$
d\mu(t) \,=\, \alpha(t)dt\,,
$
in which  $\alpha(t)$ is the  integrable function
$$
\alpha(t)= \frac{\partial m}{\partial t}\, (t,\bar{x}(t))\,.
$$
For functions $m(t,x)$'s that are merely of bounded variation with respect to $t$, the idea is to define the integral according to (\ref{partial_variation1}), but now taking $\mu(\cdot)$ to be some measure constructed from limits of  finite difference approximations of the function $m(\cdot,\cdot)$.  
The first step is to define the measure $\mu(\cdot)$ to replace 
$ \frac{\partial m}{\partial t}\, (t,\bar{x}(t))$. The construction of the measure is based on the following lemma, which invokes the hypotheses:
\begin{itemize}
\item[(BV1):] $\bar{x}(\cdot)$ is continuous and there exists 
$\delta' >0$ such that $m(t,.)$ is continuously differentiable on the interior of $\bar{x}(t)+ \delta' \,\B$ for all $t\in [S,T]$.
\item[(BV2):] 
\begin{itemize}
\item[(i):]  $m(\cdot,x)$ has bounded variation along $\bar x(\cdot)$
\item[(i):] $\nabla_{x} m(\cdot,x)$ has bounded variation along $\bar x(\cdot)$.
\end{itemize}
\end{itemize}
\begin{prop}
\label{4.1}
Take functions $m(\cdot,\cdot):[S,T] \times \R^{n} \rightarrow \R^{r}$ and $\bar x(.): [S,T] \rightarrow \R^{n}$. Assume (BV1) and (BV2) are satisfied. Take any sequence of partitions $\{ t_{i}^{j} \}_{i=0}^{N_{j}} $, $j=1,2,\ldots$ of $[S,T]$ such that diam$(\{ t^{j}_{i} \})  \rightarrow 0$ as $j \rightarrow \infty$ and any sequence $\rho^{j} \downarrow 0$. Take also any sequence of collections of $n$-vectors  $\{ \xi^{j}_{i} \}_{i=0}^{N_{j}} $  such that
$$
 \xi^{j}_{i} \,\in\, \bar x(t)+\rho^{j} \B  \quad \mbox{for some } t \in [t^{j}_{i},t^{j}_{i+1} ]
$$
for each $i$ and $j$. Define the sequence of discrete measures $\mu^{j}(\cdot)$, $j=1,2,\ldots$  to be
$$
\mu^{j}(t)\,=\, \sum_{i=0}^{N_{j}-1}\left[ 
m(t^{j}_{i+1}, \xi^{j}_{i})- m(t^{j}_{i}, \xi^{j}_{i}) 
 \right] \delta(t-t^{j}_{i})
$$
(`pseudo-density' notation), in which $\delta(\cdot)$ denotes the Dirac delta function.
\ \\

\noindent
Then there exists a (signed) Borel measure $\mu(\cdot)$ on $[S,T]$ such that
$$
\mu^{j}(\cdot) \rightarrow \mu(\cdot) \quad \mbox{with respect to the weak$^{*}$ topology on $C^{*}([S,T]; \R^{r})$, i.e.}
$$
$$
\int_{[S,T]} g^{T}(t)\,d \mu^{j}(t) \rightarrow \int_{[S,T]} g^{T}(t)\,d \mu(t)
\mbox{ for every $g(.) \in C([S,T]; \R^{r})$}\,.
$$
Furthermore, the limit measure $\mu(\cdot)$ does not depend on the choice of sequences of partitions $\{ t^{j} \}_{i=0}^{N_{j}} $, the sequence $\{\rho_{j}\}$ or the sequence of collections of  vectors $\{ \xi^{j} \}_{i=0}^{N_{j}}$ satisfying the stated conditions.
\end{prop}
\vspace{0.0 in}

\noindent
{\bf Proof.} Denote by $\eta(\cdot)$ and $\tilde{\eta}(\cdot)$ cummulative variation functions of $ m(\cdot,x)$ and $\nabla_{x} m(\cdot,x)$, respectively, along $\bar x(\cdot)$. We can choose $\bar{\epsilon} >0$ and $\bar \delta >0$ such that $\eta^{\bar \delta}_{\bar \epsilon}(T) < \infty$, $\tilde{\eta}^{\bar \delta}_{\bar \epsilon}(T) < \infty$ and
$
\theta_{\bar{x}}(\bar{\epsilon}) \,\leq\, \bar \delta\,,
$
where $\theta_{\bar x}(\cdot)$ is a continuity modulus for $\bar x(\cdot)$.
\vspace{0.05 in}

\noindent
Take  any sequence of partitions $\{ t^{j}_{i} \}$ of $[S,T]$, $j=1,2,\ldots$, any sequence of real numbers $\rho^{j}\downarrow 0$  and any sequence $\{ \xi^{j}_{i} \}$, $j=1,2,\ldots$ of collections of  $n$-vectors with the properties listed in the statement of the lemma. Write
$$
\epsilon_{j}:= \sup_{j' \geq j} \, \mbox{diam}(\{ t^{j'}_{i}\}) \,.
$$
By assumption $\epsilon_{j} \downarrow 0$, as $j \rightarrow \infty$.
Fix  $j$ and $j' (> j)$. Consider the case
\begin{itemize}
\item[(T):]  $\{ t^{j}_{i} \} \subset  \{ t^{j'}_{i} \}$\,, i.e. $\{t^{j'}_{i} \}$ is a sub-partition of $\{t^{j}_{i} \}$\,.
\end{itemize}
We relabel the sequence  $\{ t^{j}_{i} \}$ as $\{s_{0}\ldots s_{M} \}$. Then, since $\{ t^{j'}_{i} \}$ is a sub-partition of  $\{ t^{j}_{i} \}$, $\{ t^{j'}_{i} \}$ can be written
 $$
\{ t^{j'}_{i} \}= \{s_{0 \ell}  \}_{\ell=0}^{ \ell_{0}}  \cup \ldots  \cup
 \{s_{(M-1) \ell}  \}_{\ell=0}^{ \ell_{M-1}} \,.
$$
Here $s_{00}=s_{0}$, and, for $i=1,\ldots, M-1$,  $s_{i0}=s_{i}$ and $s_{i \ell_{i}}= s_{i+1}(= s_{(i+1)0})$. Relabel the $n$-vectors associated with these two partitions as $\xi_{0}, \ldots \xi_{M-1}$ and as $\xi_{i\ell}$, $\ell=0,\ell_{i}-1$, $i=1, \dots, M-1$.
For any continuous function $g(.):[S,T] \rightarrow \R^{r}$
\begin{eqnarray*}
\< \mu^{j}-\mu^{j'}, g(\cdot) \> &:=&
\int_{[S,T]}g(t)(d\mu^{j}(t)- d\mu^{j'}(t))
\\
&=& \sum_{i=0}^{M-1}\left[ 
g^{T}(s_{i})\left( m(s_{i+1},\xi_{i}) -m(s_{i},\xi_{i}) \right) \right.
\\
&& 
\hspace{0.2 in}   - \sum_{\ell=0}^{\ell_{i}-1}\left.
g^{T}(s_{i \ell})\left( m(s_{i(l+1)},\xi_{i \ell}) -m(s_{i \ell},\xi_{i \ell}) \right)
\right]\,.
\end{eqnarray*}
Using the fact that
$
|m(s_{i(\ell +1)}, \xi_{i \ell}) - m(s_{i\ell}, \xi_{i \ell})| \,\leq\, \eta^{\bar{\delta}}_{\bar{\epsilon}}(s_{i( \ell +1)}) - \eta^{\bar{\delta}}_{\bar{\epsilon}}(s_{i\ell})\,, \mbox{etc.,}
$
\vspace{0.05 in}

we can write
\begin{equation}
\label{decomp}
\< \mu^{j}-\mu^{j'}, g(\cdot) \>\,=\,a + e_{1}\,,
\end{equation}
where
\begin{eqnarray}
\nonumber
&&a = \left.\sum_{i=0}^{M-1}
g^{T}(s_{i})\right[
\left( 
m(s_{i+1},\xi_{i}) -m(s_{i},\xi_{i}) \right)
\\
&&
\label{a_term} 
\hspace{1.3 in}  
 -\left. \sum_{\ell=0}^{\ell_{i}-1} 
\left( m(s_{i(l+1)},\xi_{i \ell}) -m(s_{i \ell},\xi_{i \ell}) \right) \right]
\,.
\end{eqnarray}
and $e_{1}$ is an `error term' that satisfies
\begin{eqnarray*}
|e_{1}| \,\leq\, \theta_{g}(\epsilon_{j})\times  \sum_{i,\ell} \left( 
\eta^{\bar \delta}_{\bar \epsilon}(s_{i(\ell +1)}) -
\eta^{\bar \delta}_{\bar \epsilon}(s_{i \ell})
\right)\,=\, \theta_{g}(\epsilon_{j})\times  \eta^{\bar \delta}_{\bar \epsilon}(T)\,.
\end{eqnarray*}
Here, $\theta_{g}(\cdot)$ is a continuity modulus for $g(\cdot)$.
Observe next that $\xi_{i \ell},\,\xi_{i} \in \bar{x}(s_{i})+ \theta_{\bar x}(\epsilon_{j})\B$ for each $i$ and $\ell=0,\ldots,l_{i}$. It follows from (BV1) that, for $j$ sufficiently large, the functions $m(s_{i (\ell+1)},\cdot)$  and  $m(s_{i \ell},\cdot)$ in (\ref{a_term}) are continuously differentiable on a ball containing $\xi_{i \ell}$ and $\xi_{i}$. Employing  an exact first order Taylor expansion of these functions about $\xi_{i}$, 
we can write terms in the inner summation  on the right of (\ref{a_term}) 
\begin{eqnarray}
\nonumber
&& 
g^{T}(s_{i})
\left( 
 m(s_{i(\ell+1)},\xi_{i \ell}) -m(s_{i \ell},\xi_{i \ell}) \right)
\\
\nonumber
&& \hspace{1.0 in} =\, 
g^{T}(s_{i})
\left( 
 m(s_{i(\ell+1)},\xi_{i }) -m(s_{i \ell},\xi_{i }) \right)
\\
\label{decompose}
&&
\hspace{1.0 in}\,+\,
g^{T}(s_{i})
\left( 
 \nabla_{x} m(s_{i(\ell+1)},\tilde{\xi}_{i \ell}) - \nabla_{x} m(s_{i \ell},\tilde{\xi}_{i \ell}) \right)\cdot \left(
\xi_{il}- \xi_{i}
 \right)
\end{eqnarray}
for $n$-vectors $\tilde \xi_{i \ell} \in \mbox{co}\,\{\xi_{i}, \xi_{i \ell}\}$. We can show that, for each $i$ and $\ell=0,\dots,l_{i-1}$, and for all values of $j$ sufficently large, 
$$| \tilde \xi_{i \ell}-
\bar x(s_{i\ell})|\leq 3 \times (\theta_{\bar x}(\epsilon_{j})+ \rho_{j})  \leq \bar \delta\,.
$$
Here, $\theta_{\bar x }(.)$ is a continuity modulus for $\bar x(.)$. We note also that 
$$ |\xi_{il} -\xi_{i}| \leq 2 \times (\theta_{\bar x}(\epsilon_{j})+ \rho_{j})\,.
$$
Substituting (\ref{decompose}) into (\ref{a_term}), noting cancellation of terms and, finally, using the fact that $ \nabla_{x}m(\cdot,x)$ has bounded variation along $\bar x(\cdot)$, we arrive at 
$$
a \,=\,0 + \ldots 0 + e_{2}\,,
$$
 where  $e_{2}$ is an error term that satisfies
$
|e_{2}| \,\leq\, 2 \,(\theta_{\bar x}(\epsilon_{j})+ \rho^{j}) \,||g(.)||_{C} \,\tilde{\eta}^{\bar \delta}_{\bar \epsilon}(T)\,.
$
\ \\

\noindent
It now follows from (\ref{decomp})  and (\ref{a_term}) that
\begin{equation}
\label{estimate}
\< \mu^{j}-\mu^{j'}, g(.) \> \,\leq \, 
\theta_{g}(\epsilon_{j})\, \eta^{\bar \delta}_{\bar \epsilon}(T)
+
2 \,(\theta_{\bar x}(\epsilon_{j})+ \rho_{j}) \, ||g(.)||_{C} \, \tilde{\eta}^{\bar \delta}_{\bar \epsilon}(T)\,.
\end{equation}
Recall that (\ref{estimate}) has been proved in the case (T). Suppose that (T) is not satisfied, i.e. $\{t^{j'}_{i}\}$ is not a sub-partition of $\{t^{j}_{i}\}$. We shall show that a similar estimate is valid.  The key observation here is that, given the two partitions, we can construct a new partition $\{\tilde{t}_{i}\}$ of $[S,T]$, simply by combining all the discretization times of the two partitions. Write $\tilde \mu(\cdot)$ for the measure
$$
\tilde{\mu}(t)\,=\, \sum_{i}\left[ 
m(\tilde{t}_{i+1}, \bar{x}(\tilde{t}_{i}))- m(\tilde{t}_{i},\bar{x}(\tilde{t}_{i})) 
 \right] \delta(t-\tilde{t}_{i})\;.
$$
Applying the preceding analysis, first to $\mu^{j}(\cdot)$ and $\tilde \mu(\cdot)$ and then  to $\mu^{j'}$  and $\tilde \mu$, and noting the triangle inequality, we arrive at:
\begin{eqnarray*}
&& |\< \mu^{j}-\mu^{j'}, g(\cdot) \> |\,\leq \, |\< \mu^{j}-\tilde{\mu}, g(\cdot) \>| +
 |\< \mu^{j'}-\tilde{\mu}, g(\cdot) \>|
\\
&& \hspace{0.5 in}
\leq 
2  \left(\theta_{g}(\epsilon_{j}) \eta^{\bar \delta}_{\bar \epsilon}(T)
+
2\,(\theta_{\bar x}(\epsilon_{j})+\rho_{j} )\, ||g(.)||_{C} \, \tilde{\eta}^{\bar \delta}_{\bar \epsilon}(T)\right)\,.
\end{eqnarray*}
This relation implies
$$
\lim_{j \rightarrow \infty} \sup_{j' \geq j}  \< \mu^{j}-\mu^{j'}, g(\cdot) \>\,=\,0\,.
$$
We have shown that, for arbitrary continuous $g(\cdot)$, $\{\< \mu^{j}, g(\cdot) \>\}$ is a Cauchy sequence in $\R$. The sequence therefore has a limit.
\ \\

\noindent
In consequence of property (\ref{basic}) of functions having 
bounded variation, the measures $\{\mu^{j}(\cdot)\}$ are bounded by $\eta^{\bar \delta}_{\bar \epsilon}(T)$,  for $j$ sufficiently large. Since closed balls in $C^{*}([S,T], \R^{r})$ are weak$^{*}$ compact
there exists a Borel measure $\mu(\cdot)$ on $[S,T]$ and a subsequence $\{\mu^{j_{k}}(\cdot)\}$  of $\{\mu^{j}\}$ such that 
$$
\mu^{j_{k}}(\cdot) \rightarrow \mu(\cdot) \quad \mbox{with respect to the weak$^{*}$ topology},
$$
as  $k \rightarrow \infty$. But then, by the preceding analysis, 
$$
\lim_{j \rightarrow \infty}\, \<\mu^{j},g(\cdot) \> = \lim_{k \rightarrow \infty}\,\<\mu^{j_{k}},g(\cdot) \>= \<\mu, g(\cdot) \>\,,
$$
for any $g(.) \in C([S,T]; \R^{r})$. We have demonstrated that there exists a Borel measure $\mu(\cdot)$ such that $\mu^{j}(\cdot)$ converges to $\mu(\cdot)$ in the manner claimed (weak$^{*}$ convergence in the dual space).
\ \\

\noindent
We now prove the final assertion of the lemma (`uniqueness of the limit'). If it were not true, there would exist two sequences of Borel measures $\{\mu^{j}(\cdot)\}$ and $\{\tilde \mu^{j}(\cdot)\}$ on $[S,T]$ that converge to different limits $\mu(\cdot)$ and $\tilde{\mu}(\cdot)$ (respectively), with respect to the weak$^{*}$ topology. The fact that the limits are distinct means that there exists some $g^{*}(.) \in C([S,T];\R^{r})$ such that
\begin{equation}
\label{different}
\< \mu, g^{*}(\cdot)\> \not= \< \tilde{\mu}, g^{*}(\cdot)\>\,.
\end{equation}
Now construct a new sequence $\{ \tilde{\tilde{\mu}}^{j}(\cdot)\}$ by alternating elements in the two sequences. By the preceding analysis, there exists a Borel measure $\tilde{\tilde{\mu}}$ such that $\tilde{\tilde{\mu^{j}}} \rightarrow \tilde{\tilde{\mu}}(\cdot)$ in the weak$^{*}$ topology, as $j \rightarrow \infty$. So
$$
\lim_{j \rightarrow \infty}\< \tilde{\tilde{\mu}}^{j}, g^{*}(\cdot)\> = \< \tilde{\tilde{\mu}}(\cdot), g^{*}(\cdot)\>\,.
$$
But the sequence  $\{\< \tilde{\tilde{\mu}}^{j}, g^{*}(\cdot)\>\}$ cannot converge,  because there exist two subsequences, one with limit 
$\< \mu, g^{*}(\cdot)\> $ and the other with limit $\< \tilde{\mu}, g^{*}(\cdot)\>$,  which are distinct by  (\ref{different}). This contradiction completes the proof. $\square$
%
\begin{defn}
Take functions $m(\cdot,\cdot): [S,T] \times \R^{n} \rightarrow \R^{r}$ and $\bar x(\cdot):[S,T] \rightarrow \R^{n}$. Assume hypothesis (BV1) and (BV2) are satisfied. Then {\it the partial variation measure of $ m(\cdot,x)$ along $\bar x(\cdot)$}, written
$$
B \rightarrow \int_{B} d_{t}m(t, \bar{x}(t)) \,,
$$
is the Borel measure on $[S,T]$:
$$
\mu(\cdot) \,=\, \lim_{j} \mu^{j}(\cdot)\,,
$$
in which the limit is taken with respect to the weak$^{*}$ topology on $C^{*}([S,T];\R^{r})$. Here, $\{\mu^{j}(.)\}$ is any sequence of discrete Borel measures, each of the form
$$
\mu^{j}\,=\, \sum_{i=0}^{N_{j}-1}\left[ 
m(t^{j}_{i+1}, \xi^{j}_{i})- m(t^{j}_{i}, \xi^{j}_{i}) 
 \right] \delta(t-t^{j}_{i})\,,
$$
in which $\{t^{j}_{i}\}_{0=1}^{N_{j}}$, $j=1,2,\ldots$, is a sequence of partitions of $[S,T]$ such that diam$(\{t^{j}_{i}\}) \rightarrow 0$ as $j \rightarrow \infty$. $\{\xi^{j}_{i}\}_{i=0}^{N_{j}-1}$ is a collection of $n$-vectors such that $\xi^{j}_{i} =\bar{x}(t)$ for some $t \in [t^{j}_{i}, t^{j}_{i+1}]$.
(The definition of $B \rightarrow \int_{B} d_{t}m(t, \bar{x}(t))$ is unambiguous since, according the preceding analysis, the limiting measure $\mu(\cdot)$ is the same for all choices of sequences $\{\mu^{j}(\cdot)\}$.)
\end{defn}

\noindent
The next proposition relates the value of the partial variation measure on a subinterval $[a,b] \subset [S,T]$ and the difference in values of $m(\cdot,x)$ at $a$ and $b$. 
\begin{prop}
\label{4_3}
Take functions $m(\cdot,\cdot): [S,T] \times \R^{n} \rightarrow \R^{r}$ and $\bar x(\cdot):[S,T]\rightarrow \R^{n}$. 
 Assume hypothesis (BV1) and (BV2) are satisfied, for some $\delta' >0$. Denote by $\eta(\cdot)$ and $\tilde{\eta}(\cdot)$ the cummulative variation functions of  $m(\cdot,x)$ and  $\nabla_{x} m(\cdot,x)$, respectively.
\ \\

\noindent
Take a closed subinterval $[a,b] \subset [S,T]$, $\delta \in (0,\delta')$ such that $\eta^{\delta}(T) < \infty$ and $\tilde{\eta}^{\delta}(T) < \infty$ and $\xi \in \R^{n}$ such that $\xi =\bar{x}(t)$ for some $t \in [a,b]$.  Assume that
\begin{equation}
\label{modulus}
\theta_{\bar x}(|b-a|) \leq \delta\,,
\end{equation}
where $\theta_{\bar x}(.)$ is a continuity modulus for $\bar x(.)$.
 Then
\begin{eqnarray*}
&&
\left|
\int_{[a,b]}d_{t}m(t, \bar x(t)) - \left( m(b,\xi) - m(a, \xi)
\right)
\right|
\,\leq \,
\\
&& \hspace{0.4 in}
\theta_{\bar x}(|b-a|) \times\left( 
\tilde{\eta}^{\delta} (b) -
 \tilde{\eta}^{\delta} (a)
\right)
+
\left( 
\eta(a)- \lim_{a' \uparrow a}  \eta (a')
\right)
+
\left(
\lim_{b' \downarrow b} 
\eta(b')-  \eta(b)
\right)\,.
\end{eqnarray*}
(The second term on the right is  interpreted as $0$ if $a=S$, and the third as $0$ if $b=T$.)
\end{prop}

\noindent
{\bf Proof.}  Fix $[a,b] \subset [S,T]$ and $\delta >0$ such that $\eta^{ \delta}(T) < \infty$ and $\tilde{\eta} ^{\delta}(T) < \infty$ and (\ref{modulus}) is satisfied.  
\ \\

\noindent
Let $\{ G_{k}(.): [S,T]\rightarrow \R^{r \times r}\}$ be a sequence continuous functions such that
\begin{eqnarray*}
&&G_{k}(t)=I_{r \times r} \mbox{ for } t \in [a,b]
\\ 
&& |G_{k}(t)|= 0 \mbox{ if } t \leq a - k^{-1} \mbox{ or }  b+k^{-1}  \leq t
\\
&& 
|G_{k}(t) |\leq 1 \mbox{ if } a - k^{-1} \leq t \leq a  \mbox{ or }  b \leq t \leq b+ k^{-1}\,.
\end{eqnarray*}
Take any index value $k$ and $\epsilon >0$ sufficiently small that   $\eta^{\delta}_{\epsilon}(T) < \infty$  and $\tilde \eta^{\delta}_{\epsilon}(T) < \infty$. Let $\{t^{j}_{i}\}_{i=0}^{N_{j}}$ be a sequence of partitions of $[S,T]$ such that diam$(\{t^{j}_{i}\})\rightarrow 0$ as $j\rightarrow 0$ and such that  $\{t^{j}_{i}\}$ contains $a$ and $b$ for each $j$. Now define
$$
\mu^{j}(t)\,=\, \sum_{i=0}^{N_{j}-1}\left[ 
m(t^{j}_{i+1}, \bar x(t^{j}_{i}))- m(t^{j}_{i}, \bar x(t^{j}_{i})) 
 \right] \delta(t-t^{j}_{i})\;.
$$
By Lemma \ref{4.1}, applied component-wise,
\begin{equation}
\label{converges}
\int_{[S,T]}G_{k}(t)\,d \mu^{j}(t) \rightarrow \int_{[S,T]} G_{k}(t)\,d_{t}m(t, \bar x(t))
\mbox{ as $j \rightarrow \infty$ }\,.
\end{equation}
For each $j$, let $m^{j}_{1}$ and $m^{j}_{2}$ be the values of the index $i$ defined by  
$t^{j}_{m^{j}_{1}}=a$ and $t^{m^{j}_{2}}_{i}=b$. Then, for each $j$, 
\begin{equation}
\label{sigma123}
\int_{[S,T]} G_{k}(t)\,d \mu^{j}(t)- (m(b,\xi)-m(a,\xi)) \,=\, \Sigma_{1}+ \Sigma_{2} + \Sigma_{3}\,,
\end{equation}
in which
\begin{eqnarray*}
\Sigma_{1}&=& \sum_{i=0}^{m^{j}_{1}-1}G_{k}(t_{i})\left[ 
m(t^{j}_{i+1},\bar x(t^{j}_{i}))- m(t^{j}_{i},\bar x(t^{j}_{i})) 
 \right]\,,
\\
\Sigma_{2}&=& \sum_{i=m^{j}_{1}}^{m^{j}_{2}-1} \left[ 
m(t^{j}_{i+1},\bar x(t^{j}_{i}))- m(t^{j}_{i},\bar x(t^{j}_{i})) 
\right]
- \left[ m(b, \xi)-m(a, \xi) \right]\,,
\\
\Sigma_{3} &=& \sum_{i=m^{j}_{2}}^{N_{j}-1}
G_{k}(t_{i})\left[ 
m(t^{j}_{i+1},\bar x(t^{j}_{i}))- m(t^{j}_{i},\bar x(t^{j}_{i})) 
 \right] \,.
\end{eqnarray*}
%
Consider the term $\Sigma_{2}$. Take any $\nu \in R^{r}$. Then, using the  exact  first order Taylor expansion formula, we can show, for $j$ sufficiently large,
\begin{eqnarray*}
&& \nu^{T} \sum_{i=m^{j}_{1}}^{m^{j}_{2}-1}
 \left[ 
 m(t^{j}_{i+1},\bar x(t^{j}_{i})) -m(t^{j}_{i},\bar x(t^{j}_{i})) 
\right]  
\\
&& \qquad
\,=\, \nu^{T} 
\sum_{i=m^{j}_{1}}^{m^{j}_{2}-1}
 \left[ 
 m(t^{j}_{i+1},\xi) -m(t^{j}_{i},\xi) 
+  
 \nabla_{x} m(t^{j}_{i+1},\xi^{\nu j}_{i}) - \nabla_{x} m(t^{j}_{i},\xi^{\nu j}_{i}
) 
\right] \cdot \left( \bar x(t^{j}_{i})- \xi \right) 
\\
&& \qquad
\,\leq\, \nu^{T} 
\left( 
\left[  
 m(b,\xi)+ 0 \ldots 0  -m(a, \xi)
\right]
\right.
\\
&&
\hspace{0.9 in} +  \;
\theta_{\bar x}(|b-a|) \times \sum_{i=m^{j}_{1}}^{m^{j}_{2}-1}
 \left[ 
 \nabla_{x} m(t^{j}_{i+1},\xi^{\nu_{j}}_{i}) - \nabla_{x} m(t^{j}_{i},\xi^{\nu j}_{i}
) 
\right]  
)
\,,
\end{eqnarray*}
in which $\xi^{\nu j}_{i} \in \bar{x}(a)+ \delta \B$, for each $i$ and all $j$ sufficiently large. (We have used the fact that $G_{k}(t)\equiv I_{r \times r}$ for $t \in [a,b]$.) But  $\nabla_{x} m(\cdot, x)$ has bounded variation along $\bar x(\cdot)$.  We can therefore conclude that, for $j$ sufficiently large,
$$
|\nu^{T}\Sigma_{2}|
\,\leq\,
|\nu| \left( 
\tilde{\eta}^{\delta}_{\epsilon}(b)- \tilde{\eta}^{\delta}_{\epsilon}(a)
\right)\times \left(\theta_{\bar x}(|b-a|) \right)\,.
$$
Since $\nu$ is an arbitrary $r$-vector, 
\begin{equation}
\label{sigma2}
|\Sigma_{2}|
\,\leq\,  
(\tilde{\eta}^{\delta}_{\epsilon}(b)- \tilde{\eta}^{\delta}_{\epsilon}(a)) 
\; 
\theta_{\bar x}(|b-a|) 
\,.
\end{equation}
Now take any $\delta' \in (0, \delta)$. Since $G_{k}(.)$ satisfies $|G_{k}(.)| \leq 1$ on $[S,T] \backslash [a,b]$ and vanishes on $[S,T] \backslash [a-k^{-1}, b+k^{-1}]$, we deduce from  property (\ref{basic})  of cummulative variation functions that, for sufficiently large $j$, 
\begin{equation}
\label{sigma12}
|\Sigma_{1}| \leq \eta^{\delta'}_{\epsilon}(a)- \eta^{\delta'}_{\epsilon}(S \vee (a- k^{-1})) \mbox{ and }
|\Sigma_{3}| \leq \eta^{\delta'}_{\epsilon}(T \vee (b+ k^{-1}))
- \eta^{\delta'}_{\epsilon}(b)\,.
\end{equation}
Noting (\ref{converges}), (\ref{sigma123}), (\ref{sigma2}) and (\ref{sigma12}) and passing to the limit as $j \rightarrow \infty$ gives
\begin{eqnarray}
\nonumber
&&|\int_{[S,T]} G_{k}(t)\,d \mu(t) - (m(b,\xi)-m(a,\xi))| \,\leq\,
\left( \tilde{\eta}_{\epsilon}^{\delta}(b)- \tilde{\eta}_{\epsilon}^{\delta}(a)\right) \times \theta_{\bar x}(|b-a|)
\\
&& 
\label{final}
\hspace{0.4 in}
 + 
( \eta^{\delta'}_{\epsilon}(a)- \eta^{\delta'}_{\epsilon}(S \vee (a- 1/k)) +  
( \eta^{\delta'}_{\epsilon}(T \vee (b+ 1/k))
- \eta^{\delta'}_{\epsilon}(b) )
\end{eqnarray}
But $\delta' >0$ and $\epsilon >0$ are arbitrary, sufficiently small numbers. We may therefore pass to the limit as first $\epsilon \downarrow 0$ and second as  $\delta' \downarrow 0$, to deduce the validity of the preceding relation when $\eta^{\delta}_{\epsilon}(\cdot)$ and $\eta^{\delta'}_{\epsilon}(\cdot)$ are replaced by $\eta^{\delta}(\cdot)$ and $\eta(\cdot)$, respectively. 
\ \\

\noindent
So far $k$ has been fixed. Finally, we pass to the limit as $k \rightarrow \infty$. Since $G_{k}(t) \rightarrow I_{r \times r}\times \chi_{[a,b]}$ everywhere and the monotone function $\eta(.)$ has everywhere one-sided limits, we deduce with the help of the Dominated Convergence Theorem that
\begin{eqnarray*}
&& |\int_{[a,b]}d_{t}m(t, \bar x(t)) - (m(b,\xi)-m(a,\xi))| \,\leq\,
\\
&& \hspace{0.3 in} \theta_{\bar x}(|b-a|) \times \left(\tilde{\eta}^{\delta}(b)- \tilde{\eta}^{\delta}(a)\right)  + 
\left( \eta(a)- \lim_{a' \uparrow a}\eta(a')\right) +  
\left( \lim_{b' \downarrow b'} \eta(b)
- \eta(b) \right)\,
\end{eqnarray*}
in which $\lim_{a' \uparrow a}\eta(a'):= \eta(S)$ if $a=S$ and $\lim_{b' \downarrow b}\eta(b'):= \eta(T)$ if $b=T$. The proof is complete. $\square$
\section{An Application}
Consider a control system relating the control function  $u(\cdot)$ to an output function $y(t)$ according to 
$$
(S)\quad
\left\{
\begin{array}{l}
\dot x(t) = f(x(t),u(t)) \mbox{ a.e } t \in [S,T]\,,
\\
u(t) \in \Omega \mbox{ a.e } t \in [S,T]\,,
\\
x(0)=x_{0}\,,
\\
y(t)=g(x(t)) \mbox{ for } t \in [S,T]\,,
\end{array}
\right.
$$
the data for which is: functions $f(\cdot,\cdot):\R^{n} \times \R^{m} \rightarrow \R^{n}$ and $g(\cdot):\R^{n} \rightarrow \R$, a set $\Omega \subset \R^{m}$ and an $n$-vector $x_{0}$.
\ \\

\noindent
Let $\bar{u}(\cdot)$ be a control function that has been chosen to give a desired value to the output at time $T$, which we write
\begin{equation}
\label{cost}
J(u(\cdot))\,:=\, g(x(T;u(\cdot),x_{0}))
\end{equation}
where $t \rightarrow x(t;u(\cdot),x_{0})$ denotes the solution to the differential equation in the control system description, for a given control function $u(\cdot)$ and initial condition $x_{0}$. (Hypotheses will be imposed ensuring the existence and uniqueness  solutions.) Write $\bar x(t)= x(t;\bar u(\cdot),x_{0})$.
\ \\

\noindent
In this section we focus our attention on the following phenomenon: in control engineering it is often the case that a feedback control cannot be implemented perfectly, but only  with a time delay. This is especially evident in process control, where controlled chemical reactors  are routinely modelled with a pure delay at the input, to take account of   the finite rate of flow of fluids between reactors, etc. (See, for example, the widely studied Tennessee Eastman challenge controller design problem, in which the system equations take the form of a matrix of first order lags with pure time delay \cite{McAvoy}). The presence of a time delay complicates the controller design and so, if it is small, it is often ignored. To justify the use of idealized `delay-free' models, it then becomes necessary to carry out a sensitivity analysis, to quantify the errors in the output $J(u(\cdot))$ when small delays are introduced into the controller implementation. We need then to look at consequences of applying the control
\begin{equation}
\label{delay_control}
u^{h}(t):=
\left\{
\begin{array}{ll}
\bar u(S) & \mbox{if }  t-h < S  
\\
\bar u(t-h) & \mbox{if } S \leq t-h \leq T 
\\
\bar u(T) & \mbox{if } T < t-h \,.
\end{array}
\right.
\end{equation}
Notice we allow $h$ to be both positive (a delay) or negative (an advance).
The effect of introducing the delay on the output at time $T$ is quantified by
$$
J(u^{h}(.)):= g(x(T;u^{h}(\cdot),x_{0}))\,.
$$
Suppose that $f(\cdot,\cdot)$ is continously differentiable and globally Lipschitz continuous. If the control $\bar u(\cdot)$ is an absolutely continuous function, a routine analysis yields the information that $h \rightarrow J(u^{h}(\cdot))$ is differentiable at the origin with gradient
\begin{equation}
\label{gradient_smooth}
\frac{d}{dh}J(u^{h}(\cdot))|_{h=0}\,=\,
\int_{[S,T]}p^{T}(t) \nabla_{u} f(\bar x(t), \bar u(t))\, \frac{d \bar u}{dt}(t)dt\, ,
\end{equation}
in which $p(.): [S,T]\rightarrow \R^{n}$ is the solution to the costate equation:
\begin{equation}
\label{adjoint}
\left\{
\begin{array}{l}
-\dot{p}(t)=  \nabla_{x}f ^{T}(\bar x(t), \bar u(t))\,p(t)
\\
p(T)= \nabla_{x}g^{T}(\bar x(T))\,.
\end{array}
\right.
\end{equation}
It is sometimes required to consider controls $\bar u(\cdot)$ that are not absolutely continuous (`bang-bang' controls arising from the solution to minimum time problems, for example). {\it Is it possible to establish regularity  properties of $h \rightarrow J(u^{h}(\cdot))$ and to derive a formula akin to (\ref{gradient_smooth}) for a larger class of controls $\bar u(\cdot)$, and when $f(x,u)$ is no longer assumed to be differentiable w.r.t. the $u$ variable}? The Prop. \ref{5.1} below provides a positive answer, when $\bar u(\cdot)$ is a function of bounded variation.  
\ \\

\noindent
We shall invoke the following hypotheses: there exists $k_{1} >0$, $\delta >0$ and a modulus of continuity $\theta(.)$
such that 
\begin{itemize}
\item[(S1):] $g(\cdot)$ is a $C^{1}$ function,
\item[(S2):]  $f(\cdot,\cdot)$ is continuous, $f(\cdot,u)$ is a $C^{1}$ function for each $u \in \Omega$ and 
\begin{itemize}
\item[(i):]
$|f(x,u)|\,\leq\, c[1 + |x|]$ for all $x \in \R^{n}$, $u \in \Omega$ 
%
\item[(ii):] $|\nabla_{x} f(x, u)| \leq K$ for all $x\in \R^{n}$ and $u \in \Omega$, 
\item[(iii):]
%
$|f(x,u)-f(x',u)-\nabla_{x} f(x',u) | \,\leq\, \theta(|x-x'|)\times |x-x'|$

\hspace{ 2.0 in}for all $x,x'  \in \bar x(t) + \delta \B$ and $u\in \Omega$,
\end{itemize}
\item[(S3):] 
%
$|f(x,u)-f(x,u')|+ |\nabla_{x} f(x,u)-\nabla_{x}f(x,u')| \leq k_{1}|u-u'|$ 

\noindent
\hspace{1.5 in} for all $x\in \bar x(t)+ \delta \B$, $u,u' \in \Omega$ and $t \in [S,T]\,.$ 
\item[(BV):] $\bar u(\cdot)$ has bounded variation.
\end{itemize}

\begin{prop}
\label{5.1}
Consider the control system (S) and a control function $\bar u(\cdot)$. Assume that
hypotheses (S1)-(S3) and (BV) are satisfied.
\ \\

\noindent
For any number $h \in \R$ define $u^{h}(\cdot): [S,T] \rightarrow \R^{m}$ according to  (\ref{delay_control}). Write $x^{h}(\cdot)$ for the solution on $[S,T]$ of $\dot x(t)= f(x(t),u^{h}(t))$, $x(S)=x_{0}$
and also
\begin{equation}
\label{mfunction}
m^{h}(t,x):= f(x, u^{h}(t))\mbox{ and } \nabla_{x} m^{h}(t,x):= \nabla_{x} f(x, u^{h}(t))\,.
\end{equation}
Then, for all $h$ in some neighborhood of $0$:
\begin{itemize}
\item[(a):]
 $m^{h}(\cdot,x)$ and $ \nabla_{x} m^{h}(\cdot,x)$ have bounded variation along $\bar x(\cdot)$, 
\item[(b):]
$h' \rightarrow J(u^{h'}(\cdot))$ (given by (\ref{cost})) has one sided derivatives (from left and right) at $h'=h$:
\begin{equation}
\label{sensitivity_1}
\lim_{h' \downarrow h}\, \frac{J(u^{h'}(.))- J(u^{h}(.))}{{h'-h}}  =  -
\int_{[S,T)} p_{h}^{T}(t)\,d_{t}m^{h}(t, \bar x(t))
\end{equation}
and
\begin{equation}
\label{sensitivity_2}
\lim_{h' \uparrow h}\, \frac{J(u^{h}(.))- J(u^{h}(.))}{{h-h'}}  = -
\int_{(S,T]} p^{T}_{h}(t)\,d_{t}m^{h}(t, \bar x(t))\,.
\end{equation}
(In these relations, $p_{h}(\cdot)$ is the solution to (\ref{adjoint}) when $u^{h}(\cdot)$ and $x^{h}(\cdot)$ replace $\bar u(\cdot)$ and $\bar x(\cdot.)$, and $B \rightarrow \int_{B}d_{t}m^{h}(t, \bar x(t))$ is the partial variation measure associated with $m^{h}(\cdot,\cdot)$.)
\item[(c):]
If $\bar u(\cdot)$ is continuous at both endpoints $S$ and $T$, the mapping $h' \rightarrow J(u^{h'}(\cdot))$ is  differentiable at  $h$ and its  derivative is
\begin{equation}
\label{sensitivity_3}
\lim_{h' \rightarrow h}\, \frac{J(u^{h'}(\cdot))- J(u^{h}(\cdot))}{{h'-h}}  = -
\int_{[S,T]} p_{h}^{T}(t)\,d_{t}m^{h}(t, \bar x(t))\,.
\end{equation}
\end{itemize}
\end{prop}

\noindent
{\bf Discussion:} 
The property that the sensitivity function $h \rightarrow J(u^{h})$ is differentiable when $\bar u(\cdot)$ has bounded variation and continuous at the two end-times (part (c) of the proposition) is highly non-trivial, since $f(x,u)$ is not assumed to be differentiable w.r.t. $u$. To convey the nature of this property in its simplest terms, let us consider the  case of control system (S) when $f(x,u)$ is independent of $x$ (write the function $f(u)$).  Assume that
\begin{itemize}
\item[(a):]
$f(\cdot)$ is Lipschitz continuous,
\item[(b):]
$\bar u(\cdot)$ is continuously differentiable.
\end{itemize}
It is straightforward to show that, under these hypotheses,   the sensitivity function $V(h):= J(u^{h}(\cdot))$ is Lipschitz continuous. A standard analysis based on perturbing $h$ and using the properties of Clarke's generalized directional derivative (c.f.  \cite[Proof of Thm. 2.7.3]{Clarkebook}) and a nonsmooth chain rule permits one to derive the following estimate of the subdifferential of $V(.)$ at $0$:
\begin{equation}
\label{sensitivity_nonsmooth}
\partial V(0) \,\subset \, - \int_{S}^{T} \,\mbox{co}\,  \partial_{u} H(p(t),\bar u(t)) \dot{\bar u} (t) dt\,.  
\end{equation}
in which 
$$
H(p,u):= p^{T}f(u)\,.
$$ 
Here $\partial_{u}H(u,p)$ denotes the subdifferential w.r.t. the $u$ variable, for fixed $p$.  (We refer to the end of Section 1 for definition the subdifferential.) $p(\cdot)$ is the solution of the adjoint system (\ref{adjoint}). 
The right side of this relation is a set valued integral, defined in the usual way as the collection of integrals of selectors of the set valued integrand.
\ \\

\noindent
Prop. \ref{5.1} tells us, contrary to what the standard analysis leading to the formula (\ref{sensitivity_nonsmooth}) might lead us to expect, the sensitivity function is actually differentiable on a neighborhood of $0$. {\it Indeed it tells us even more: the sensitivity function is differentiable on a neighborhood of $0$, when $f(t,x)$ is $x$-dependent and $u(.)$ is merely a function of bounded variation (continuous at its endpoints). This  surprising regularity property, is a consequence of the properties of functions of bounded variation, along a specified trajectory, established in earlier sections of this paper.}
%
\ \\

\noindent
{\bf Proof of Prop 5.1.} 
For any $h$ we have:
\begin{eqnarray*}
&&\int_{[S,T]} |f(\bar x(t), u^{h}(t)) -f(\bar{x}(t), \bar{u}(t))|dt
\\
&& \hspace{0.2 in} \leq k_{1}\left(\int_{[(S+ h)\wedge T,T]}|u(t)- u(t-h)| +
\int_{[S,(S+ h)\wedge T]}|u(t)- u((S))|
\right)
\\
&& \hspace{0.2 in}
=k_{1}\left( \int_{[(S+ h)\wedge T,T]}\eta_{\bar u}(t)dt+ \int_{[S, (S+h)\wedge T]}\eta_{\bar u}(t)dt \right)
\leq 2k_{1}\; \eta_{\bar u}(T) \; h\,.
\end{eqnarray*}
in which $\eta_{\bar u}(\cdot)$ is the cummulative variation function of $u(\cdot)$. By Filippov's Existence Theorem (see e.g. \cite[Thm. 2.4.3]{Vinter}), there exist a number $K_{1}$, independent of $h$, such that
\begin{equation}
\label{filippov}
||x^{h}(.)- \bar x(.)||_{L^{\infty}} \leq K_{1} \; h\,. 
\end{equation}
 Under hypothesis (S2) there exists a unique solution to the differential equation $x^{h}(.)$ for each $h$.
  It can be deduced from (S3) that there exists $\bar h>0$ such that   $m^{h}(\cdot,x)$ 
 and $\nabla_{x}m^{h}(\cdot,x)$
  have bounded variation along $x^{h}(\cdot)$,  for all $h\in [-\bar h, \bar h]$. Write the cummulative variation functions $\eta^{h}(\cdot)$ and $\tilde \eta^{h}(\cdot)$, and their $\delta$-perturbed versions $\eta^{ h \delta}(\cdot)$ and $‪\eta^{ h \delta}(\cdot)$ respectively.
\ \\

\noindent 
We now examine the one-sided differentiability properties of $h' \rightarrow J(u^{h'})$ at $h=0$ and derive  the given formula for the one-sided derivative from the left.
(Analogous arguments can be used to treat the other cases.)
\ \\

 \noindent 
Take an arbitrary sequence $h_{i} \downarrow 0$. Then, for each $i$,
\begin{eqnarray*}
&& J(u^{h_{i}}(\cdot))- J(\bar u(\cdot)) 
\\
&& \hspace{0.3 in}=\; g(x^{h_{i}}(T))- g(\bar x(T)) - \int_{[S,T]}p^{T}(t) \left[(\dot{x}^{h_{i}}(t) - \dot{\bar{x}}(t)) \right.
\\
&& \hspace{1.8 in} 
\left. - \left( f(x^{h_{i}}(t), u^{h_{i}}(t)) -f(\bar{x}(t), \bar{u}(t) \right) \right] dt\,.
\end{eqnarray*}
A routine analysis, in which we make use of the costate equation  and right boundary condition (\ref{adjoint})  on $p(\cdot)$, apply integration by parts to the integral $\int_{[S,T]} p^{T}(t)(\dot{x}^{h_{i}}(t) - \dot{\bar{x}}(t))dt$ and consider first order Taylor expansions of $g(\cdot)$ about $\bar x(T)$ and of $x \rightarrow f(x, u^{h_{i}}(t))$ about $\bar x(t)$, reveals that
\begin{eqnarray}
\label{A1}
&& h_{i}^{-1}\left( J(u^{h_{i}}(\cdot))- J(\bar u(\cdot)) \right)
\\
\nonumber
&& \hspace{0.3 in}=\;  h_{i}^{-1}  \int_{[S,T]}p^{T}(t) \left[  f(\bar x(t), u^{h_{i}}(t)) -f(\bar{x}(t), \bar{u}(t) \right] dt 
+ e(h_{i})
\\
\nonumber
&& \hspace{0.3 in}=\;  h_{i}^{-1} \int_{[S,T]}p^{T}(t) \left[ m^{0}((t-h_{i}) \vee S), \bar x(t)) - m^{0}(t, \bar x(t)) \right] dt + e(h_{i})\,,
\end{eqnarray}
in which the `error term'  $e(h_{i})$ satisfies $|e(h_{i})| \leq ||p(.)||_{L^{\infty}}\; (\theta (K_{1} h_{i}) \; K_{1}+ \theta_{g}(K_{1}h)K_{1})$. ($\theta(.)$ is the continuity modulus of (S2) and $\theta_{g}(.)$ is a `second order' continuity modulus for $g(\cdot)$, i.e. a function such  that $|g(x)-g(x') -\nabla g(x') | \leq \theta_{g}(|x-x'|) |x-x'|$ for all $x,x' \in \bar x (T)+ \delta \B$.) We see
\begin{equation}
\label{B1} 
e(h_{i}) \rightarrow 0 \mbox{ as }i \rightarrow \infty\,.
\end{equation}
Take any  $\delta >0$ such that $\tilde \eta^{\delta}(T) < \infty$. 
Since $\eta^{\delta}(.)$ is a monotone function, there exists a subset   ${\cal I} \subset (S,T)$, of full Lebesgue measure, on which $ \eta^{\delta}(.)$ is continuous. 
%
   Prop. \ref{4_3} tells us that, for each $t \in {\cal I}$ and $i$ sufficiently large, 
\begin{eqnarray}
\nonumber
&&m^{0}((t-h_{i})\vee S),\bar x(.)) - m(t, \bar x(t)) 
\\
 &&
\label{C1}
\qquad
=\; - \int_{[(t-h_{i})\vee S, t]}d_{s}m^{0}(s, \bar x(s)) + e_{1}(t,h_{i})
\end{eqnarray}
in which the `error term' $e_{1}(t,h_{i})$ satisfies
\begin{eqnarray}
\label{D1}
 |e_{1}(t, h_{i})| &\leq & \theta_{\bar x}(|h_{i}|)\left(  \tilde \eta^{\delta}(t)- \tilde \eta^{\delta}((t-h_{i})\vee S)  \right)\,.
\end{eqnarray}
It follows that
\begin{eqnarray}
 \nonumber
&& h_{i}^{-1}\int_{[S,T]}|p^{T}(t) 
 e_{2}(t, h_{i})|dt 
\\
\nonumber
&& \quad\leq\,
h_{i}^{-1}\theta_{\bar x}(|h_{i}|)\; ||p(.)||_{L^{\infty}} \;
\int_{[S,T]} \left( 
\tilde{\eta}^{\delta} (t) -
 \tilde{\eta}^{\delta} ((t-h_{i}) \vee S)
\right)dt
\\
\label{E1}
&& \quad=\; \theta_{\bar x}(|h_{i}|) \; ||p(.)||_{L^{\infty}} 
\; \left( h^{-1}_{i}\int_{[T-h_{i},T]}\tilde \eta^{\delta}(t)dt + \tilde \eta(S)  \right) \; \rightarrow\;  0 \,,
\end{eqnarray}
as $i \rightarrow \infty$. In consequence of Fubini's Theorem
\begin{eqnarray*}
&&
h^{-1}_{i}\int_{[S,T]} \int_{[(t-h_{i})\vee S, t]}p^{T}(t)d_{s}m^{0}(s, \bar x(s))dt
\,=\, \int_{[S,T]}p_{i}^{T}(s) d_{s}m^{0}(s, \bar x(s))ds\,,
\end{eqnarray*}
in 
which 
$$
p_{i}(s):= h_{i}^{-1}\int_{[s,(s+h_{i})\wedge T]}p(t)dt\quad \mbox{for } s \in [S,T]\,.
$$
Since $p(.)$ is continuous, $p_{i}(t)\rightarrow \tilde p(t)$ for all $t \in [S,T]$, where
\begin{equation*}
\tilde p(t)\,\rightarrow \,
\left\{
\begin{array}{ll}
p(t)& \mbox{if } t \in [S,T)
\\
0 & \mbox{if } t= T\;.
\end{array}
\right.
\end{equation*}
By the Dominated Convergence Theorem
\begin{eqnarray}
\nonumber
\lefteqn{h^{-1}_{i}\int_{[S,T]} \int_{[(t-h_{i})\vee S, t]}p^{T}(t)d_{s}m^{0}(s, \bar x(s))dt}
\\
\nonumber
& \hspace{0.5 in} \rightarrow& \, \int_{[S,T]}\tilde p^{T}(s) d_{s}m^{0}(s, \bar x(s))ds
\\
\label{F1} 
& \hspace{0.5 in} =& \, \int_{[S,T)}p^{T}(s) d_{s}m^{0}(s, \bar x(s))ds\,.
\end{eqnarray}
Combining relations (\ref{A1}) - (\ref{E1})), we arrive at
$$
h_{i}^{-1}\,\left (J(u^{h_{i}}(\cdot))- J(\bar u(\cdot)) \right) \,\rightarrow\, -
\int_{[S,T)}p^{T}(s) d_{s}m^{0}(s, \bar x(s))ds\,.
$$
We have confirmed formula (\ref{sensitivity_1}) and the existence of the limit. 
\ \\

\noindent
We now attend to the final assertion of the proposition. Suppose then that $\bar u(\cdot)$ is continuous at $S$ and $T$. For $h$ sufficiently small,  $m^{h}(., x)$ is continuous at $t=S$ and $t=T$, uniformly as $x$ ranges over neighborhoods of $x^{h}(S)$ and $x^{h}(T)$ . It can be deduced from relation (\ref{z1}) in Prop. \ref{prop2.3} 
that $\eta(.)$ is right and left continuous at $S$ and $T$. Then, by Prop. \ref{4_3}, $d_{t}m^{h}(\cdot,x^{h}(t))$ has no atom at either $S$ or $T$. The differentiability of $h' \rightarrow J(u^{h'}(\cdot))$  and the formula (\ref{sensitivity_3}) now follow from (\ref{sensitivity_1}) and (\ref{sensitivity_2}), since the integrals in the latter two formulae, over $[S,T)$ and $(S,T]$ respectively,  are the same.  $\square$
\newpage

\begin{center}
{\Large {\bf Appendix: Proof of Prop. \ref{prop2.3}}}
\end{center}
Take   $\bar \epsilon >0$ such that $\eta^{\bar \delta}_{\bar \epsilon}(T) < \infty$ and $\delta \in (0, \bar \delta)$. To begin, we verify the following assertion: for any $\delta' \in (0, \delta)$
\begin{equation}
\label{right_cont}
 \eta^{\delta'}(T) \,\leq\, \underset{T' \uparrow T}{\lim}\, \eta^{\delta}(T')
+ \underset{x \in \bar x(T) + \delta \B, a \in A}{\sup}
d_{H}(F(T,x,a), F(T^{-},x,a))
,
\end{equation}
%
Take any $\delta'' \in (\delta', \delta)$. Since $\bar x(.)$ is continuous, we can choose $T_{1} \in (S,T)$ such that 
\begin{eqnarray}
\label{inclA}
\bar x([T_{1},T])+ \delta' \B &\subset& \bar x(T)+  \delta'' \B,
\\
\label{inclB}
\bar x([T_{1},T])+ \delta'' \B &\subset& \bar x(T)+ \delta \B
\end{eqnarray}
 
For any $T_{2} \in (T_{1}, T]$ define 
$$
G_{\delta''}(T_{1},T_{2})\,:=\,
\sup\,\left\{ 
\sum_{i=0}^{N} \underset{x \in \bar x(T) + \delta'' \B, a \in A}{\sup}
d_{H}( F(t_{i+1},x,a), F(t_{i},x,a))
\right\}\,,
$$
in which the `outer' supremum is taken over all partitions $\{ t_{i}\} $ of $[T_{1}, T_{2}]$. Notice that the `inner' suprema are all taken over the same set $(\bar x(T) + \delta \B)\times A$. It follows that the value of $G_{\delta''}(T_{1},T_{2})$ is unchanged if we restrict the magnitude of the diameters of the partitions considered in the definition; that is, for any $\epsilon >0$ we have
\begin{equation}
\label{restrict_eps}
G_{\delta''}(T_{1},T_{2})\,=\,
\sup\,\left\{ 
\sum_{i=0}^{N} \underset{x \in \bar x(T) + \delta'' \B, a \in A}{\sup}
d_{H}( F(t_{i+1},x,a),  F(t_{i},x,a))\,|\,\mbox{ diam}\, (\{t_{i}\})\, \leq \epsilon
\right\}\,.
\end{equation}
%
By considering the modification of arbitrary partitions of $[S,T]$ to include the extra `grid point' $T_{1}$ and taking account of (\ref{inclA}) and (\ref{measure}), we see that, for any $\epsilon \in (0, \bar \epsilon ]$,
$$
\eta_{\epsilon}^{\delta'} (T)\,\leq\, 
\eta_{\epsilon}^{\delta'}(T_{1})+ 2 \, \gamma (c  \,\epsilon) + G_{\delta''}(T_{1},T)\,.
$$
Here, $c$ and $\gamma(.)$ are as in hypotheses (C1) and (C2). In the limit as $\epsilon \downarrow 0$, we obtain
\begin{equation}
\label{first}
\eta^{\delta'}(T)\,\leq\,\eta^{\delta'}(T_{1})+ G_{\delta''}(T_{1},T)\,.
\end{equation}
Take any  $\rho > 0$.  Then there exists a partition $\{ t^{\rho}_{0},\ldots, t^{\rho}_{N^{\rho}} \}$ of $[T_{1}, T]$ achieving the `outer' supremum defining $G_{\delta''}(T_{1},T)$, with error at most $\rho$. From (\ref{inclB}), we deduce
\begin{equation}
\label{second}
 G_{\delta''}(T_{1},T)\,\leq\,  G_{\delta''}(T_{1},t^{\rho}_{N^{\rho}-1})+ 
\underset{x \in \bar x(T) + \delta'' \B,\, a \in A}{\sup}
d_{H}(F(t^{\rho}_{N^{\rho}-1},x,a), F(T,x,a)) + \rho\,.
\end{equation}
Because $\{t^{\rho}_{i}\}$ is a partition, $t^{\rho}_{N^{\rho}-1} < T$. In view of (\ref{restrict_eps}) which, we recall, is valid for any $\epsilon>0$, we can arrange that 
\begin{equation}
\label{close}
|T- t^{\rho}_{N^{\rho}-1}| \leq \rho\,.
\end{equation}
But, in consequence of (\ref{inclB}) and (\ref{restrict_eps}), we know that, for any $\epsilon \in (0, \bar \epsilon]$, 
$$
\eta_{\epsilon}^{\delta} (t^{\rho}_{N^{\rho}-1})\,\geq\, 
\eta_{\epsilon}^{\delta}(T_{1})+ G_{\delta''}(T_{1},t^{\rho}_{N^{\rho}-1})\,.
$$
In the limit, as $\epsilon \downarrow 0$, we obtain
\begin{equation}
\label{third}
\eta^{\delta} (t^{\rho}_{N^{\rho}-1})\,\geq\, 
\eta^{\delta}(T_{1})+ G_{\delta''}(T_{1},t^{\rho}_{N^{\rho}-1})\,.
\end{equation}
Combining  (\ref{first}),(\ref{second}) and (\ref{third}) and noting that 
$\eta^{\delta'}(T_{1}) \leq \eta^{\delta}(T_{1})$, we see that
\begin{eqnarray*}
\eta^{\delta'}(T)&\leq&  \eta^{\delta}(t^{\rho}_{N^{\rho}-1})+
\underset{x \in \bar x(T) + \delta'' \B,\, a \in A}{\sup}
d_{H}(F(t^{\rho}_{N^{\rho}-1},x,a), F(T,x,a)) 
+ \rho\,.
\end{eqnarray*}
This relation is valid for any $\rho >0$.  Passing to the limit as $\rho \downarrow 0$, while taking account of  (\ref{close}) and using Prop. \ref{prop1} (part (b)) to evaluate the limit of the sup term on the right side, and noting that $\delta'' < \delta$, we deduce
$$
\eta^{\delta'} (T)\,\leq\,  \lim_{T' \uparrow T} \eta^{\delta}(T')+ \underset{x \in \bar x(T) + \delta \B, a \in A}{\sup}
d_{H}(F(T^{-},x,a), F(T,x,a))\,.
$$
This confirms relation (\ref{right_cont}).
%
%
\ \\

\noindent
%
Our next task will be to relate the $\delta$-perturbed cummulative variation functions $\eta^{\delta}(\cdot)$ and $\tilde{\eta}^{\delta}(\cdot)$ of 
$ F(\cdot,x,a)$ and $\tilde F(\cdot,x,a)$, respectively, at times $t \in (S,T)$.  
\ \\

\noindent
Fix $t \in (S,T)$
and take $\epsilon \in (0, \bar{\epsilon}]$. 
Let ${\cal T}= \{t_{0}=S,\ldots, t_{N}=t\}$ be an arbitrary partition of $[S,t]$ with diam$({\cal T})\leq \epsilon$. Take an arbitrary sequence $s_{j}\downarrow S$. Then, for $j$ sufficiently large, $s_{j}< t_{1} $ and 
\begin{eqnarray*}
&&\eta^{\delta}_{\epsilon}(t)\,\geq\, \sup \left\{ d_{H}(F(s_{j},x,a) ,F(S,x,a) )\;|\; x \in \bar x([S,s_{j}])+\delta B, a\in A \right\}
\\
&&\hspace{1.0 in} + \sup \left\{ d_{H}(F(s_{j},x,a) ,F(t_{1},x,a) )\;|\; x \in \bar x([s_{j},t_{1}])+\delta B, \, a\in A \right\}
\\
&&\hspace{1.0 in} + \sum_{i=1}^{N-1}\sup \left\{ d_{H}(F(t_{i+1},x,a) ,F(t_{i},x,a) )\;|\; x \in \bar x([t_{i},t_{i+1}])+\delta B, a\in A \right\}\;.
\end{eqnarray*}
In view of Prop. \ref{prop1}, we may pass to the limit as $j\rightarrow \infty$ in this relation to obtain:
\begin{eqnarray*}
&&\eta^{\delta}_{\epsilon}(t)\,\geq\, \sup \left\{ d_{H}(F(S^{+},x,a) ,F(S,x,a) )\;|\; x \in \bar x(S)+\delta B,\, a\in A\right\}
\\
&&\hspace{1.0 in} + \sup \left\{ d_{H}(F(S^{+},x,a) ,F(t_{1},x,a) )\;|\; x \in \bar x([S,t_{1}])+\delta B, \, a\in A \right\}
\\
&&\hspace{1.0 in} + \sum_{i=1}^{N-1}\sup \left\{ d_{H}(F(t_{i+1},x,a) ,F(t_{i},x,a) )\;|\; x \in \bar x([t_{i},t_{i+1}])+\delta B, \, a\in A \right\}\;.
\end{eqnarray*}
Since ${\cal T}$ was an arbitrary  partition with  diam$({\cal T}) \leq \epsilon$, it follows that 
\begin{equation}
\label{one}
\eta^{\delta}_{\epsilon}(t) \geq \sup \left\{ d_{H}(F(S^{+},x,a) ,F(S,x,a) )\;|\; x \in \bar x(S)+\delta B, \, a\in A\right\} + \tilde \eta^{\delta}_{\epsilon}(t)\,.
\end{equation}
Take any partition ${\cal T}= \{t_{0}=S, \ldots, t_{N}=t\}$ of $[S,t]$ with  diam$({\cal T}) \leq \epsilon$. Then
\begin{eqnarray}
\nonumber
&&\tilde \eta^{\delta}_{\epsilon}(t)\,\geq\, \sup \left\{ d_{H}(F(S^{+},x,a) ,F(t_{1},x,a) )\;|\; x \in \bar x([S,t_{1}])+\delta B, a\in A \right\}\qquad \qquad
\\
\label{three*}
&& + \sum_{i=1}^{N-1}\sup \left\{ d_{H}(F(t_{i+1},x,a) ,F(t_{i},x,a) )\;|\; x \in \bar x([t_{i},t_{i+1}])+\delta B, a\in A \right\}.
\end{eqnarray}
By the triangle inequality we have, for each $x\in \bar x([S,t_{1}]) +\delta \B$ and $a\in A$,
$$
d_{H}(F(S^{+},x,a) ,F(t_{1},x,a) )\,\geq\, d_{H}(F(S,x,a) ,F(t_{1},x,a) ) - d_{H}(F(S^{+},x,a) ,F(S,x,a))\;.
$$
Furthermore,
\begin{eqnarray*}
&& \sup d_{H}(F(S^{+},x,a) ,F(t_{1},x,a) )\geq
\\
&& \hspace{0.5 in}  \sup d_{H}(F(S,x,a) ,F(t_{1},x,a) ) - \sup\{ d_{H}(F(S^{+},x,a) ,F(S,x,a) )\,,
\end{eqnarray*}
where, in each term, the sup is taken over $(x,a)\in (\bar x( [S,t_{1}])+ \delta \B) \times A$.
Since ${\cal T}$ was an arbitrary partition such that diam$( {\cal T} )\leq \epsilon$,  we deduce from (\ref{three*})
that
\begin{eqnarray}
\nonumber
&&\tilde \eta^{\delta}_{\epsilon}(t)\,\geq\, \eta^{\delta}_{\epsilon}(t) - \sup \{ d_{H}(F(S^{+},x,a) ,F(S,x,a) )\,|\,(x,a)\in (\bar x([S,t_{1}]) +\delta \B) \times A    \}\,.
\end{eqnarray}
This relation combines with (\ref{one}) to yield
\begin{eqnarray}
\nonumber0&\leq&
\eta^{\delta}_{\epsilon}(t)-\tilde \eta^{\delta}_{\epsilon}(t) - \sup\{ d_{H}(F(S^{+},x,a) ,F(S,x,a) )\,|\,x\in \bar x(S)+ \delta \B, \, a\in A\}
\\
\label{Delta}
&\leq& \Delta(\epsilon, \delta)\,,
\end{eqnarray}
in which
\begin{eqnarray}
\nonumber
&&  
\hspace{-0.2 in} \Delta(\epsilon, \delta)\,:=\,
\sup\{ d_{H}(F(S^{+},x,a) ,F(S,x,a) )\,|\,(x,a) \in (\bar x([S,(S+\epsilon)\wedge T]) +\delta \B) \times A
\\ 
&&
\label{Delta_e_delta}
\hspace{0.3 in}
-\,\sup \{ d_{H}(F(S^{+},x,a) ,F(S,x,a) )\,|\,x\in \bar x(S)+ \delta \B, a\in A\}\;.
\end{eqnarray}
Since $\bar x(.)$ is continuous and, as is easily shown, $F(S^{+},\cdot,\cdot)$ has modulus of continuity $\gamma(\cdot)$ on $\{\bar{x}(S)+\delta \B\} \times A$, where $\gamma(\cdot)$ is as in hypothesis (C2), we have
$$
\lim_{\epsilon' \downarrow 0} \Delta(\epsilon', \delta)=0\;.
$$
Considering (\ref{Delta}) in the limit as $\epsilon \downarrow 0$ yields
\begin{equation}
\label{additional}
\eta^{\delta}(t)=\tilde \eta^{\delta}(t) + \max \{ d_{H}(F(S^{+},x,a) ,F(S,x,a) )\,|\,x\in \bar x(S) + \delta \B, \, a\in A\}\;.
\end{equation}
Taking the limit $\delta \downarrow 0$ yields 
\begin{equation}
\label{1.1}
\eta(t)=\tilde \eta(t) +\underset{a\in A}{\sup} \,  d_{H}(F(S^{+},\bar{x}(S),a) ,F(S,\bar{x}(S),a) )\;.
\end{equation}
We have validated the property  (\ref{z2}) in the case $t \in (S,T)$.
\ \\

\noindent
Now we verify the right continuity of $\tilde{\eta}^{\delta}(\cdot)$  and $\tilde \eta(\cdot)$ at $S$. We shall merely confirm the right continuity of  $\tilde \eta^{\delta}(\cdot)$  at $S$, since  this will imply the right continuity of $\tilde \eta(\cdot)$ at $S$, in consequence of the relation
$$
0 \leq\ \tilde \eta(t)-\tilde \eta(S) =\tilde \eta(t)-0\leq \tilde \eta^{\delta}(t)-0= \tilde \eta^{\delta}(t)-
 \tilde \eta^{\delta}(S)\; \mbox{ for all } t \in (S,T].
$$
Let us assume, in contradiction, that $\tilde \eta^{\delta}(.)$ is not right continuous at $S$. Then there exists $\alpha > 0$ such that
$
\tilde \eta^{\delta}(t)-(\tilde \eta^{\delta}(S)=0) \geq \alpha \mbox{ for all $t \in (S,T]$}
$. Taking account of Prop. \ref{prop1}, we can choose  $\epsilon \in (0,T-S)$ such that $\tilde \eta^{\delta}_{\epsilon}(t) < \infty$ and
\begin{equation}
\label{bound1}
 \sup_{(x,a) \in (\bar{x}([S,t])+ \delta B)\times A} d_{H}(F(t,x ,a),F(S^{+},x,a ) ) \leq \alpha /4
\end{equation}
for all $t \in [S,S+\epsilon]$.
We may now choose a partition $\{s_{0}, ..., s_{N}\}$  of $[S,S+\epsilon]$, of diameter at most $\epsilon$,  such that
\begin{eqnarray}
\nonumber
\tilde \eta^{\delta}_{\epsilon}(S+\epsilon) &\leq&  \sup_{(x,a) \in  (\bar{x}([S,s_{1}])+ \delta B)\times A }d_{H}(F(s_{1},x ,a),F(S^{+},x,a ) ) 
+\Sigma_{2} + \alpha /4
\\
\label{inequality}
&\leq&   \alpha/4 + \Sigma_{2}+\alpha/4 \,=\, \Sigma_{2} + \alpha/2\,, 
\end{eqnarray}
in which 
$$
\Sigma_{2}\,:=\, \sum_{i=1}^{N-1}\sup_{(x,a) \in (\bar{x}([s_{i},s_{i+1}])+ \delta \B) \times A} d_{H}(F(s_{i+1},x,a ),F(s_{i},x,a ) ) \,.
$$
Now we choose a partition $\{t_{0}, \ldots,t_{M}\}$ of $[S,s_{1}]$, with diam$(\{t_{i}\})\leq \epsilon$,   such that
$$
\alpha \leq \tilde \eta^{\delta}_{\epsilon}(s_{1}) \leq \Sigma_{1}+ \alpha/4 \,,
$$
where
$$
\Sigma_{1}\,:=\, \sum_{i=0}^{N-1}\sup_{\underset{a\in A}{x \in \bar{x}([t_{i},t_{i+1}])+ \delta \B}} d_{H}(F(t_{i+1},x ,a),F(t_{i},x,a ) ) \,.
$$
It follows that 
$$
\Sigma_{1} \geq 3\alpha/4\,.
$$ 
But since the concatenation of $\{t_{0}\ldots,t_{M}\}$ and $\{s_{1}\ldots,s_{N}\}$ is a partition of   $[S,S+\epsilon]$, of diameter no greater than $\epsilon$, we know from the preceding inequality that
$$
\tilde{\eta}^{\delta}_{\epsilon}(S+ \epsilon) \geq \Sigma_{1} + \Sigma_{2} \geq  \Sigma_{2} + 3\alpha/4\,.
$$
But this contradicts (\ref{inequality}). We have confirmed that $\tilde \eta^{\delta}(.)$ (and so also $\tilde \eta(.)$) are continuous from the left at $S$.
\ \\

\noindent
Next, we shall show that
\begin{eqnarray}
\label{right_cont1}
&& \eta^{\delta}(T) \,=\, \underset{T' \uparrow T}{\lim}\, \eta^{\delta}(T')
+ \underset{x \in \bar x(T) + \delta \B, a \in A}{\sup}
d_{H}(F(T,x,a), F(T^{-},x,a))\,,
\\
&&
\label{right_cont2}
 \eta(T) \,=\, \underset{T' \uparrow T}{\lim}\, \eta(T')
+ \underset{a \in A}{\sup}\;
d_{H}(F(T,\bar x(T),a), F(T^{-},\bar x(T),a))\,.
\end{eqnarray}
This will complete the proof of the remaining assertions of the proposition. Indeed, since the multifunction $t \rightarrow \tilde F(t,.,.)$  is continuous at the right end-point T, the analysis leading to 
(\ref{additional}) and (\ref{1.1})
, but now applied to $\tilde F$, yields
$$
\tilde \eta^{\delta}(T) \,=\, \underset{T' \uparrow T}{\lim}\,\tilde  \eta^{\delta}(T')
+ 0\quad  \mbox{ and }   \quad \tilde \eta(T) \,=\, \underset{T' \uparrow T}{\lim}\, \tilde \eta(T')
+ 0\,.
$$
This is the claimed left continuity of $\tilde \eta^{\delta} (\cdot)$ (and, by implication,  of $\tilde \eta(.)$) at $T$. On the other hand, (\ref{right_cont1}) and (\ref{right_cont2}) combine with (\ref{additional}) and (\ref{1.1}) to yield the representation of $\tilde \eta^{\delta}(T)$  in terms of  $\eta^{\delta}(T)$  in (\ref{z2}) (and the analogous representation of $\tilde \eta(T)$).
\ \\

\noindent
To prove  (\ref{right_cont1}) and (\ref{right_cont2}) we first note that, since 
%
$\delta' \rightarrow \eta^{\delta'}(T)$ is monotone, we can find $\delta_{1} \in (\delta, \bar \delta)$, arbitrarily close to $\delta$, such that $\delta' \rightarrow \eta^{\delta'}(T)$ is continuous at $\delta_{1}$. 
But then,  by (\ref{right_cont}),
\begin{eqnarray}
\nonumber
 \eta^{\delta_{1}}(T)&=&\underset {\delta' \uparrow \delta_{1}}{\lim} \eta^{\delta_{1}}(T) \,
\\
\label{right_cont3}
&\leq& \underset{T' \uparrow T}{\lim}\, \eta^{\delta_{1}}(T')
+ \underset{x \in \bar x(T) + \delta_{1} \B, a \in A}{\sup}
d_{H}(F(T,x,a), F(T^{-},x,a))\,.
\end{eqnarray}
By Lemma \ref{simple} however we have, for any $T' \in (S,T)$,  
\begin{equation}
\label{basic_bound} 
\eta^{\delta}(T)- \eta^{\delta}(T') \,\leq\,
\eta^{\delta_{1}}(T)- \eta^{\delta_{1}}(T').
\end{equation}
Consequently
$$
\eta^{\delta}(T)- \lim_{T'\uparrow T} \eta^{\delta}(T') \,\leq\,
\eta^{\delta_{1}}(T)-  \lim_{T'\uparrow T}\eta^{\delta_{1}}(T').
$$
This relation combines with (\ref{right_cont3}) to give
\begin{equation}
\label{right_cont5}
 \eta^{\delta}(T) \,\leq\, \underset{T' \uparrow T}{\lim}\, \eta^{\delta}(T')
+ \underset{x \in \bar x(T) + \delta_{1} \B, a \in A}{\sup}
d_{H}(F(T,x,a), F(T^{-},x,a)) \,.
\end{equation} 
Since $\delta_{1}$ can be chosen such that  $\delta_{1}-\delta$ is arbitrarily small and in view of the continuity properties of $F(.,.,.)$, we see that the preceding relation is true when the supremum is taken over $x \in \bar x(T)+\delta \B$ in place of $x \in \bar x(T)+\delta_{1} \B$.
\ \\

\noindent
Take $\epsilon \in (0, T-S)$ and $\epsilon' \in(0, \epsilon)$. Since,  $\eta_{\epsilon'}^ \delta(T)\leq \eta_{\epsilon}^ \delta(T)$, we have
\begin{eqnarray}
\label{limitaa}
\eta_{\epsilon}^ \delta(T)&\geq & 
\eta_{\epsilon'}^{\delta}(T-\epsilon)
+ \underset{x \in \bar x([T-\epsilon, T]) + \delta \B, a \in A}{\sup}
d_{H}(F(T,x,a), F(T-\epsilon,x,a)) \,.
\end{eqnarray}
Passing to the limit in the preceding inequality, first as $\epsilon' \downarrow 0$ and then as  $\epsilon \downarrow 0$, noting that
\begin{eqnarray*}
&& \underset{\epsilon \downarrow 0}{\lim}\; \left\{\underset{x \in \bar x([T-\epsilon, T]) + \delta \B, a \in A}{\sup}
d_{H}(F(T,x,a), F(T-\epsilon,x,a))\right\}
\\
&& \hspace{2.0 in}  =\;
 \underset{x \in \bar x(T) + \delta \B, a \in A}{\sup}
d_{H}(F(T,x,a), F(T^{-},x,a))\, ,
\end{eqnarray*}
and taking account of  (\ref{right_cont5}) (when $\delta$ replaces $\delta_{1}$), we deduce
(\ref{right_cont1}).
\ \\

\noindent
Now take $\delta' \in (0, \delta)$. From (\ref{limitaa})
\begin{eqnarray*}
\eta_{\epsilon}^ \delta(T)&\geq & 
\eta_{\epsilon'}^{\delta'}(T-\epsilon)
+ \underset{x \in \bar x([T-\epsilon, T]) + \delta \B, a \in A}{\sup}
d_{H}(F(T,x,a), F(T-\epsilon,x,a)) \,.
\end{eqnarray*}
Passing to the limit, as $\epsilon' \downarrow 0$, $\delta' \downarrow 0$ and $\epsilon \downarrow 0$ and  $\delta \downarrow 0$ (in that order), yields
\begin{eqnarray}
\label{limitaaa}
\eta(T)&\geq & 
\underset{T' \uparrow T}{\lim}   \eta(T')
+ \underset{ a \in A}{\sup}\;
d_{H}(F(T,\bar x(T),a), F(T^{-},\bar x(T),a)) \,.
\end{eqnarray}
Take $\delta_{1} \in (0, \bar \delta)$ such that $\delta \rightarrow \eta^{\delta}(T)$ is continuous at $\delta_{1}$.
Passing to the limit as $\delta \downarrow 0$ in (\ref{basic_bound}), for fixed $T'$, and then as $T' \rightarrow T$  yields
\begin{equation*}
\eta(T)- \lim_{T' \uparrow T} \eta (T') \,\leq\,
\eta^{\delta_{1}}(T)- \lim_{T\ \uparrow T} \eta^{\delta_{1}}(T').
\end{equation*}
But then from (\ref{right_cont3}) we deduce
\begin{equation*}
 \eta(T) \,\leq\, \underset{T' \uparrow T}{\lim}\, \eta(T')
+ \underset{x \in \bar x(T) + \delta_{1} \B, a \in A}{\sup}
d_{H}(F(T,x,a), F(T^{-},x,a)) \,.
\end{equation*}
Since  $\delta_{1}$, in this relation can be chosen arbitrarily small, the relation remains valid when we set $\delta_{1}=0$. Taking note also of (\ref{limitaaa}), we conclude (\ref{right_cont2}). The proof is complete. $\square$
\ \\

\noindent
\noindent
{\bf Acknowledgments.} This work was co-funded by
the European Union under the 7th Framework Programme
``FP7-PEOPLE-2010-ITN'', grant agreement number
264735-SADCO, and by EPSRC  under the grant “Control For Energy and Sustainability”, grant reference EP/G066477/1.

\end{document}